\documentclass[a4paper,twoside,12pt,leqno]{article}
\usepackage{amsmath,amsthm,amssymb,enumerate}

\usepackage{eucal}

\usepackage{epsfig, float, afterpage, array}

\swapnumbers
\theoremstyle{plain}
\newtheorem{theorem}{Theorem}[section]
\newtheorem{lemma}[theorem]{Lemma}
\newtheorem{corollary}[theorem]{Corollary}
\newtheorem{proposition}[theorem]{Proposition}

\theoremstyle{definition}
\newtheorem{definition}[theorem]{Definition}
\newtheorem{exercise}[theorem]{Exercise}
\newtheorem{definitions}[theorem]{Definitions}
\newtheorem{example}[theorem]{Example}
\newtheorem{examples}[theorem]{Examples}
\newtheorem{notation}[theorem]{Notation}

\newtheorem{remark}[theorem]{Remark}
\newtheorem{remarks}[theorem]{Remarks}

\newtheorem{hypotheses}[theorem]{Hypotheses}

\newtheorem{observations}[theorem]{Observations}
\newtheorem{assumptions}[theorem]{Assumptions}
\newtheorem{review}[theorem]{Review}
\numberwithin{equation}{theorem}

\makeatletter
\newcommand\Lsub{\@startsection{equation}{0}{0pt}{-\baselineskip}{0\baselineskip}{\raggedright\textbf}}
\makeatother

\newcommand{\abs}[1]{\left\lvert#1\right\rvert} %absolute value
\newcommand{\gen}[1]{\langle#1\rangle}

\def \bargamma{\Gamma}
\def \hbar {\Omega}
\DeclareMathOperator{\rank}{rank}
\DeclareMathOperator{\barr}{ \bar r}
\DeclareMathOperator{\PSL}{PSL}
\DeclareMathOperator{\height}{\alpha_3}%inf 3
\DeclareMathOperator{\depth}{\beta_2}%sup 2

\DeclareMathOperator{\fun}{\theta}
\DeclareMathOperator{\blocks}{blocks}

\DeclareMathOperator{\upp}{\sigma}

\DeclareMathOperator{\rep}{rep}
\DeclareMathOperator{\Sym}{Sym}
\DeclareMathOperator{\Alt}{Alt}
\DeclareMathOperator{\core}{core}
\def \naturals{\mathbb{N}}
\def \integers {\mathbb{Z}}

\def\into{\hookrightarrow}
\def\onto{\twoheadrightarrow}

\def\d1{\discretionary{-}{}{-}}

\begin{document}

\pagestyle{myheadings} \markboth{On the intersection of free
subgroups in free products of groups} {Warren Dicks and S.\ V.\
Ivanov}

\title{On the intersection of free subgroups \\ in free products of groups}

\author
{Warren Dicks and S.\ V.\ Ivanov}

\date{\small\today}

\maketitle

 \hfill{\it Dedicated to the memory of Prof. Charles Thomas}
 \bigskip

\begin{abstract}
Let $(G_i \mid i \in I)$ be a family of groups, let $F$ be a free group,
and let $G = F \ast \mathop{\text{\Large $*$}}\limits_{i\in I} G_i$,
the free product of $F$ and all the $G_i$.

Let $\mathcal{F}$ denote the set of all  finitely generated
subgroups  $H$ of $G$  which have the property that, for each $g \in G$ and each $i \in I$,
$H \cap G_i^{g} =  \{1\}$.
By the Kurosh Subgroup Theorem, every element of $\mathcal{F}$ is a free group.
For each free group~$H$,   the {\it reduced rank} of $H$, denoted $\barr(H)$, is defined as
 $\max \{\rank(H) -1, 0\}\,\,\in\,\,\naturals \cup \{\infty\} \subseteq [0,\infty].$
To avoid the vacuous case, we make the additional assumption
that $\mathcal{F}$ contains a non-cyclic group, and we define
\vskip -0.6cm\null \begin{center}
$\upp\,\,  :=\,\, \sup \{\textstyle\frac{\barr(H\cap K)} {\barr (H) {\cdot}\barr(K)}
: H, K \in \mathcal{F} \text{ and } \barr (H) {\cdot} \barr(K) \ne 0\}\,\,\, \in \,\,\,[1,\infty]$.\end{center}
\vskip -0.6cm\null

\noindent We are interested in precise bounds for $\upp $.
In the special case where $I$ is empty,
Hanna Neumann proved that $\upp \in [1,2]$, and
conjectured that $\upp = 1$; almost fifty years later, this interval has not been reduced.

With the understanding that $\frac{\infty}{ \infty -2}$ is~$1$, we define
\vskip -0.6cm\null \begin{center}
$\fun\,\, := \,\,\max\{ \textstyle \frac{\abs{L}}{ \abs{L}-2}
 : \text{$L$ is a subgroup of $G$ and }
\abs{L} \ne 2 \}\,\,\, \in \,\,\,[1,3].$\end{center}
\vskip -0.6cm\null

 Generalizing Hanna Neumann's theorem, we prove that
$\upp \in [\fun,\,2\fun]$, and, moreover, $\upp  = 2\fun$
whenever $G$ has 2-torsion. Since $\upp$ is finite,
$\mathcal{F}$ is closed under finite intersections.
Generalizing Hanna Neumann's conjecture, we conjecture that
$\upp  =  \fun$ whenever $G$ does not have 2-torsion.

\bigskip

{\footnotesize
\noindent \emph{{\normalfont 2000}\,Mathematics Subject Classification.} Primary: 20E06;
Secondary: \!20F06, 20F32.

}

\end{abstract}

\section{Outline}
\label{sec:intro}

\hskip.7cm  Let us first record the conventions and notation
that we shall be using.

 Throughout the article, let $G$ be a group.
 Except where otherwise specified, our $G$-actions will be on the left.

\begin{definitions}\label{defs:begin} To indicate disjoint unions, we shall use the
 symbols $\vee$, $\bigvee$ in place of
$\cup$, $\bigcup$.

We let $\naturals$
denote the set of finite cardinals, $\{0,1,2,\ldots\}$.
For each set $S$, we define  $\abs{S} \in \naturals \vee \{\infty\} \subseteq [0,\infty]$
to be the cardinal of~$S$ if $S$ is finite, and to be $\infty$ if $S$ is infinite.

For any $n \in   \{1,2,3,\ldots\}\vee\{\infty\}$, we let
$C_n$ denote a multiplicative cyclic group of order~$n$.
For any $n \in   \naturals$, we let  $\Sym_n$ denote the
group of permutations of $\{1,2,\ldots,n\}$, and we let $\Alt_n$
denote the subgroup of even permutations.

Let $a$, $b$ be elements of $G$, and let $S$ be a subset of
$G$. We shall denote the inverse of $a$ by $\overline a$. Also,
$b^a:= \overline a b a$, $\overline S := \{\overline c\mid c
\in S\}$, and  $S^a = \{c^a \mid c \in S\}$.

The {\it rank} of $G$ is defined as
$$\rank(G):= \min\{\abs{S} : \text{$S$ is a generating set of $G$} \}
\in \naturals \vee \{\infty\} \subseteq [0,\infty].$$ If $G$ is
a free group, the {\it reduced rank} of  $G$ is defined as
$$\barr(G):= \max\{\rank(G) -1, 0\} \in \naturals \vee
\{\infty\} \subseteq [0,\infty];$$ thus, $\barr(G) \,\, =\,\,
b_1^{(2)}(G)$, the first $L^2$-Betti number of  $G$; see, for
example,\linebreak \cite[Example~7.19]{Lueck02}.

\smallskip

Define
 $\height(G):= \,\,\,\inf\{ \abs{L} : \text{$L$
is a subgroup of $G$ and $\abs{L} \ge 3$}\};$ it is understood
that the infimum of the empty set is  $\infty$. By the Sylow
Theorems, $\height(G)$ is $\infty$ or $4$ or an odd prime.

 Let $\fun$ denote the
bijective, strictly decreasing (or orientation-reversing) function
$\fun\colon[3,\infty] \to [1,3], \,\,\, x\mapsto\textstyle\frac{x}{x-2}.$
Let $\fun\!\height(G)\,\,:= \,\, \fun(\height(G))= \textstyle\frac{\height(G)}{\height(G)-2}$; thus,
$$\fun\!\height(G) \,\,\in \,\,\{
\textstyle\frac{3}{1}, \textstyle\frac{4}{2}, \textstyle\frac{5}{3}, \frac{7}{5}, \frac{11}{9}, \ldots, \textstyle\frac{\infty}{\infty-2}\}\,\, = \,\,
\{1, \ldots, \frac{11}{9}, \frac{7}{5}, \frac{5}{3},  2, 3\} \,\,\subseteq \,\,[1,3].$$
For example:  $\fun\!\height(G)= 3$ if $G$ has a subgroup of order $3$;
$\fun\!\height(G) = \frac{7}{5}$ if $G$ has a subgroup of order $7$
 but none of order $3$, $4$, or $5$; and $\fun\!\height(G) = 1$ if every finite subgroup of $G$ has
 order at most $2$.
It is easy to see that if $\abs{G}
\ge 3$, then $\fun\!\height(G)\,\, = \,\,\max\{ \textstyle
\frac{\abs{L}}{\abs{L}-2} : \text{$L$ is a subgroup of $G$ and
} \abs{L} \ne 2 \}.$

Finally, define $\depth(G) :=
\begin{cases}
2 &\text{if $G$ has a subgroup of order two,}\\
{{}} 1 &\text{otherwise.}
\end{cases}$

One could define $\depth(G)$ as $\sup\{  \abs{L} : \text{$L$ is a subgroup of $G$ and }
\abs{L} \le 2 \}$ to mirror the definition of $\height(G)$.\hfill\qed
\end{definitions}

Our main interest in this article is the following.

\begin{notation}\label{not:main}
Let $(G_i \mid i \in I)$ be a family of groups, let $F$ be a free group,
and let $G = F \ast \mathop{\text{\Large $*$}}\limits_{i\in I} G_i$,
the free product of $F$ and all the $G_i$.

For each $j \in I$, we write
$G_{\neg j}:= F \ast \mathop{\text{\Large $*$}}\limits_{i\in I -\{j\}} G_i$, which gives
 $G = G_j \ast G_{\neg j}$.

Let $\mathcal{F}$ denote the set of all  finitely generated
subgroups  $H$ of $G$ which have the property that,
for each $g \in G$ and each $i \in I$,
$H \cap G_i^{g} =  \{1\}$. It follows from Kurosh's classic
Subgroup Theorem~\cite[Theorem~I.7.8]{DicksDunwoody89}
 that every element of $\mathcal{F}$ is a free group; see,
for example,~\cite[Theorem~I.7.7]{DicksDunwoody89}.

To avoid the vacuous case, we assume that some element of $\mathcal F$ has rank at least two.
We then define
\begin{equation}\label{eq:kappa}
\upp(\mathcal F)\quad=\quad  \sup \{\textstyle\frac{\barr(H\cap K)} {\barr (H){\cdot}\barr(K)}
\mid H, K \in \mathcal{F}, \, \, \barr (H){\cdot}\barr(K) \ne 0\}\quad \in\quad [1,\infty];
\end{equation}
notice that $\upp(\mathcal F) \ge 1$ since $\mathcal F$ contains some free group $H$ of rank two,
and, then, for $K = H$, we have $\textstyle\frac{\barr(H\cap K)} {\barr (H){\cdot}\barr(K)} = \frac{1}{1{\cdot}1}$.
\hfill\qed
\end{notation}

\begin{observations}  Suppose that Notation~\ref{not:main} holds.

We are interested in bounds for $\upp(\mathcal F)$.

\vskip -1cm\null
\Lsub{}\label{rem:sec11a} \hskip -0.25cm \textbf{Remarks.}  Consider the case where $I$ is empty.

Here, $G$ is a free group, $\mathcal{F}$ is the set of all finitely generated (free) subgroups of $G$,
and $\depth(G) = \fun\!\height(G) = 1$.

Let us write $\upp = \upp(\mathcal F)$.

In 1954, in~\cite{Howson}, A. G.\ Howson proved that $\sigma
\in [1,5]$, and, hence, the intersection of any two finitely
generated subgroups of a free group is again finitely
generated, that is, $\mathcal{F}$ is closed under finite
intersections. In 1956, in~\cite{HN1},   Hanna Neumann proved
that $\upp  \in [1,3]$; then, in~1958, in~\cite{HN2}, she
proved  that~$\upp  \in [1,2]$ and she conjectured that $\upp =
1$.  Almost fifty years later, the interval has not been
reduced any further, although the conjecture has received much
attention; see, for
example,~\cite{Burns71},~\cite{Dicks94},~\cite{DF01},~\cite{Imrich77},~\cite{Gersten83},~\cite{Stallings83},~\cite{Tardos92},~\cite{Tardos96}.

We now return to the general case.

\vskip -1cm\null
\Lsub{}\label{rem:sec11b} \hskip -0.25cm \textbf{Remarks.}
Let us write $\upp = \upp(\mathcal F)$, $\beta  = \depth(G)$ and $\fun = \fun\!\height(G)$.

We conjecture that $\upp  = \beta{\cdot}\!\fun$.

In Theorem~\ref{th:6},
 we prove that $\upp  \in [\beta{\cdot}\!\fun\!,\,\,2 {\cdot}\!\fun]$.

In the case where $G$ has 2-torsion, that is, $\beta=2$, then
$\upp  = \beta{\cdot}\!\fun = 2{\cdot}\!\fun$, and this case of the conjecture is true.

In the case where $G$ is 2-torsion free, that is $\beta = 1$,
then $\upp  \in [\fun,\,\,2 {\cdot}\!\fun];$
this generalizes Hanna Neumann's Theorem.  Here, our conjecture reduces to
$ \upp = \fun$, which generalizes Hanna Neumann's Conjecture.

Since $2{\cdot}\!\fun$ is finite, $\mathcal{F}$ is closed under finite intersections.
This generalizes Howson's Theorem.  An even more general statement can be deduced
from the proof of~\cite[Theorem~2.13(1)]{Sykiotis05}; see Remarks~\ref{rems:concl}(iv), below.
See also~\cite[Theorem~2]{Ivanov99} for the case
where $F$ is trivial.

\vskip -1cm\null
\Lsub{}\label{rem:sec121} \hskip -0.25cm \textbf{Remarks.}
 The condition that some element of $\mathcal F$ has rank at least two implies the following.

For each $j \in I$, $\abs{G_{\neg j}} \ge 2$.

Moreover, if,
for some  $j \in I$, $\abs{G_{\neg j}} = 2$, then
there exists a unique $j' \in I -\{j\}$ such that $\abs{G_{j'}} = 2$
and, here,  $\abs{G_{\neg j'}} \ge 3$.

\vskip -1cm\null
\Lsub{}\label{rem:sec12} \hskip -0.25cm \textbf{Remark.}
 The condition that some element of $\mathcal F$ has rank at least two is
equivalent to the condition that exactly one of the following holds.
\begin{enumerate}[(i).]
\vskip-0.75cm \null
  \item All the $G_i$ are trivial and $\rank(F) \ge 2$.
\vskip-0.8cm \null
  \item  There exists some $i_0 \in I$ such that $\abs{G_{i_0}} \ge 2$ and
$\abs{G_{\neg {i_0}}} \ge 3$.
\end{enumerate}

\vskip -1.3cm\null
\Lsub{}\label{rem:sec13} \hskip -0.25cm \textbf{Remarks.}
By the Kurosh Subgroup Theorem, again,
 each finite subgroup of $G$ lies in a conjugate of some $G_i$;
see, for example,~\cite[Proposition~I.7.11]{DicksDunwoody89}.
Hence, if $I$ is nonempty, then
\newline \null\hskip 1cm $\height(G) = \min\{\height(G_i) \mid i \in I\} \quad \text{and} \quad
\fun\!\height(G) = \max\{\fun\!\height(G_i) \mid i \in I\};$ \newline
we can arrange for $I$ to be nonempty by adding a trivial group to the family.

\vskip -1cm\null
\Lsub{}\label{rem:sec14} \hskip -0.25cm \textbf{Remark.}
 In the case where each $G_i$ is a torsion group, $\mathcal{F}$ is the
set of {\it all} finitely generated free subgroups of $G$.
\hfill\qed
\end{observations}

The organization of the paper is as follows.

In Section~\ref{sec:lower}, we  use Euler characteristics and
Bass--Serre theory, see  \cite{Bass93}, \cite{Serre80},  \cite{DicksDunwoody89}, to
show that $\upp(\mathcal{F}) \ge
\depth(G){\cdot}\fun\!\height(G)$.

%%StoW I inserted ", see  \cite{Bass93}, \cite{Serre80}," in the sentence above.
%%StoW  Perhaps, reference to  \cite{DicksDunwoody89}
%%StoW  could also be added here.
%%WtoS  Ok.  Done.  
%%WtoS  When I refereed \cite{Bass93}, Bass wrote to the editor that it
%%WtoS  was obvious that the referee was Tits or Serre. :) 
%%
%%WtoS  You can delete all this now.

  Let  $A$ and $B$ be finite subsets of $G$ with at least two elements each.
By a  {\it single-quotient} subset of $A \times B$, we mean any subset  $C$ with the property that
$\abs{\{a \overline b \mid (a,b) \in C\}} = 1$.
Sections~\ref{sec:3},~\ref{sec:4}, and~\ref{sec:5} are devoted to proving
 Corollary~\ref{cor:inequalities}(ii)
which says that, if  $\mathcal{C}$ is a set of pairwise-disjoint, single-quotient
subsets of $A \times B$, then
$\textstyle\sum\limits_{C \in \mathcal{C}} (\abs{C} - 2) \,\, \le \,\,
\fun\!\height(G){\cdot}(\abs{A} -2){\cdot}(\abs{B} -2).$

In Section~\ref{sec:freeproduct}, we use the latter result and
Bass--Serre theory
to show that $\upp(\mathcal{F}) \le
2{\cdot}\!\fun\!\height(G)$. As in the extension of Hanna
Neumann's theorem by W.~D.~Neu\-mann~\cite{WN}, we find that
all the results remain valid if, in the definition
of~$\upp(\mathcal F)$ in~\eqref{eq:kappa}, we replace
$\barr(H\cap K)$ with $\sum\limits_{s \in S} \barr(H^s\cap K)$
for any set $S$ of
 $(H,K)$-double coset representatives in $G$; see Theorem~\ref{th:5}, below.

\section{Lower bounds}\label{sec:lower}

In this section, in Proposition~\ref{prop:lower}, we prove that, if Notation~\ref{not:main}
 holds, then $\upp(\mathcal{F}) \ge \depth(G){\cdot}\!\fun\!\height(G)$.

The following is standard; see, for example,~\cite[Definition~IV.1.10]{DicksDunwoody89}.

\begin{review}\label{rev}  Suppose that $G$ is (isomorphic to) the fundamental
group of a finite graph of finite groups, $\pi(G(-),Y,Y_0)$.

 We write $VY$ and $EY$ for the vertex-set and edge-set of $Y$,
respectively.

The {\it Euler characteristic} of $G$ is defined as
$$\textstyle{\chi(G) = (\sum\limits_{v \in VY} \frac{1}{\abs{G(v)}}) -
(\sum\limits_{e \in EY} \frac{1}{\abs{G(e)}}).}$$

By Bass--Serre Theory, if $L$ is any subgroup of $G$ of finite
index, then $L$ is also the fundamental group of some finite
graph of finite groups, and\linebreak $\chi(L) =
(G:L){\cdot}\chi(G)$.

There exists a normal subgroup $H$ of $G$ of finite index
such that, for each $v \in VY$, the composite  $G(v) \into G \onto G/H$
is injective.  Moreover, any such subgroup $H$ is a finitely generated free group,
and $\chi(H) = 1 - \rank(H)$.
Thus, if $\chi(G) < 0$, then $0 > (G:H){\cdot}\chi(G)  = \chi(H) = -\barr(H)$.
\hfill \qed
\end{review}

For the purposes of this section, we introduce the following.

\begin{notation}
If $G$ contains a free subgroup of rank~2,
we let $\upp(G)$ denote the value given by $\sigma(\mathcal{F})$ in~\eqref{eq:kappa}
when $\mathcal{F}$ is taken to be the set of all finitely generated free subgroups of $G$.
\hfill\qed
\end{notation}

\begin{proposition}\label{prop:chars}  Suppose that $G$ is the fundamental
group of a finite graph of finite groups and that $\chi(G) <
0$.
If $H$ and $K$ are free normal subgroups of $G$ of finite index
such that $HK = G$, then $\barr(H \cap K) =
\frac{-1}{\chi(G)}{\cdot}\barr(H){\cdot}\barr(K) > 0$, and,
hence, $\upp(G) \ge \frac{-1}{\chi(G)}$.
\end{proposition}

\begin{proof}  Notice that $(G:K) = (HK:K) = (H:H\cap K)$, since $H\cap K$ is the
 kernel of the induced map $H \onto HK/K$. Hence,
\begin{align*}
&\chi(H){\cdot}\chi(K) =(G:H) {\cdot}\chi(G){\cdot}(G:K){\cdot}\chi(G) =
(G:H)  {\cdot}\chi(G){\cdot}(H:H\cap K) {\cdot}\chi(G) \\
&= (G:H\cap K) {\cdot}\chi(G){\cdot}\chi(G)
= \chi(H \cap K) {\cdot}\chi(G).
\end{align*}
Since $\chi(G) < 0$, we have
 $$ (-\barr(H)){\cdot}(-\barr(K)) =  \chi(H){\cdot}\chi(K) =  (-\barr(H \cap K)){\cdot}{\chi(G)} > 0 ,$$
and the result follows.  The hypothesis that $H$ is a normal
subgroup can be omitted.
\end{proof}

We now consider four concrete examples which will be used in the proof of Proposition~\ref{prop:lower}.

\begin{example}\label{ex:222} Let   $G = C_2 \ast C_2 \ast C_2$.

Then
$\chi(G) = \frac{1}{2} + \frac{1}{2} + \frac{1}{2} - 1- 1 = \frac{-1}{2}$.

We have a presentation $G = \gen{x,y, z \mid x^2 = y^2 = z^2 = 1}$.

In $\Sym_2$, consider $x' = y' = z' = (1,2)$.  There is an induced homomorphism
$G \to \Sym_2$ which sends $w$ to $w'$ for each $w \in \{x,y,z\}$.  Let $H$ be the kernel of
this homomorphism.  As in Review~\ref{rev}, $H$ is a free normal subgroup of $G$ of finite
index.  Notice that $H$ contains $xy$ and $xz$.

In $\Sym_4$, consider $x'' = (1,2)$, $y'' = (3,4)$, $z'' = (1,2)(3,4)$.
There is an induced homomorphism
$G \to \Sym_4$ which sends $w$ to $w''$ for each $w \in \{x,y,z\}$.  Let
$K$ be the kernel of this homomorphism.   As in Review~\ref{rev}, $K$ is
a free normal subgroup of $G$ of finite index.  Notice that $K$ contains $xyz$.

Then $HK$ contains $xy$, $xz$ and $xyz$.  It follows that $HK = G$.
By Proposition~\ref{prop:chars}, $\upp(G) \ge \frac{-1}{\chi(G)} = 2$.
This was also shown in~\cite[Theorem 3]{Ivanov99}.
\hfill\qed
\end{example}

In the three remaining examples, we shall tacitly use analogous constructions of free normal subgroups
of~$G$ of finite index, $H$ and $K$.

\begin{example}\label{ex:2V} Let   $G = C_2 \ast V$ where $V = C_2 \times C_2$.

Then  $\chi(G) = \frac{1}{2} + \frac{1}{4} - 1 = \frac{-1}{4}$.

We have a presentation $G = \gen{x,y,z \mid x^2 = y^2 = z^2 = (yz)^{2} = 1}$.

In $\Sym_4$, consider $x' = (1,2)(3,4)$, $y' = (1,2)$, and $z'=(3,4)$.
The resulting kernel $H$ contains $xyz$, $(xy)^2$ and $(xz)^2$.

In  $\Sym_4$, consider $x'' = (1,3)$, $y'' = (1,2)$, and $z'' = (3,4)$.
Here,
 $$x''y'' = (1,2,3)\quad \text{and} \quad x''z'' = (1,3,4).$$
The resulting kernel $K$  contains $(xy)^3$ and $(xz)^3$.

Then $HK$ contains $xyz$, $(xy)^2$, $(xz)^2$, $(xy)^3$ and $(xz)^3$.  It follows that\linebreak $HK = G$.
By Proposition~\ref{prop:chars}, $\upp(G) \ge \frac{-1}{\chi(G)} = 4$.
\hfill\qed
\end{example}

\begin{example}\label{ex:2p} Let $p$ be $4$ or an odd prime, and let $G = C_2 \ast C_p$.

Then
$\chi(G) = \frac{1}{2} + \frac{1}{p} - 1 = \frac{2-p}{2p}$.

We have a presentation $G = \gen{x,y \mid x^2 = y^p = 1}$.

Let $q = \begin{cases}2 &\text{ if $p = 4$,}\\
p &\text{ if $p$ is an odd prime.}
\end{cases}$

In $\Sym_{q+2}$, consider  $x' = \begin{cases}(1,3)(2,4) &\text{ if $p = 4$,}\\
(p+1,p+2) &\text{ if $p$ is an odd prime,}
\end{cases}$ \newline
\null\hskip 4.4cm $y' = (1,2,\ldots,p-1,p)$.
\newline Then $x'y' = y'x'$.  The resulting kernel $H$ contains $(xy)^{2p}$ and
$x(xy)^q$.

In $\Sym_{2p}$, consider
$x'' = (1,p+1)(2,3),$ \newline
\null\hskip 4.1cm $y'' = (1,2,\ldots,p)(p+1,p+2,\ldots,2p)$. \newline
Then $x''y'' = (1, 3, 4,\ldots,  p, p+1, \ldots, 2p-1, 2p).$
The resulting kernel  $K$ con\-tains~$(xy)^{2p-1}$.

Then,  $HK$ contains $(xy)^{2p}$, $x(xy)^q$ and $(xy)^{2p-1}$.
 It follows that  $HK = G$.
By Proposition~\ref{prop:chars},  $\upp(G) \ge \frac{-1}{\chi(G)} = \frac{2p}{p-2}$.
For $p=3$,  this was also shown in\linebreak~\cite[Theorem 1]{Ivanov01}.

\vskip -0.95cm\null
\Lsub{} \hskip -0.25cm \textbf{Remarks.}
For $p\ge 4$, the foregoing $K$ has rather large rank.

For $p= 4$, an alternative $K$ can be constructed by taking, in $\Sym_4$,\linebreak
 $x'' = (1,2)$, $y'' = (1,2,3,4)$.
Then, $x''y'' = (2,3,4)$ and $K$ contains $(xy)^3$.
Here, $3$ is coprime to $2p$.

For $p \ge 5$, an alternative $K$ can be constructed by taking, in  $\Sym_{p+1}$, \linebreak
$x'' \colon t \mapsto -\frac{1}{t}$, $y'' \colon t \mapsto t+1$, where we
identify $\{1,\ldots,p+1\}$ with the projective line
over the field with $p$ elements, $\mathbb{F}_p\vee \{\infty\}$.
Then, $x''y'' \colon t \mapsto -\frac{1}{t+1}$ and $K$ contains $(xy)^3$.
Here, $3$ is coprime to $2p$.

\null\vskip -1.47cm\null
\Lsub{} \hskip -0.25cm \textbf{Remarks.}
For $p = 3$, there are interesting examples
related to the action of the arithmetic group
$\PSL_2(\integers) \simeq C_2 \ast C_3$  by M\"{o}bius transformations
on the upper half-plane
$\mathfrak{h}$, the set of
complex numbers with positive imaginary part.

For $n \in \naturals$, let  ${\bargamma}_n$ denote the kernel of the mod-$n$  map
$\PSL_2(\integers) \onto \PSL_2(\integers_n)$.
For $n \ge 2$,
$ \bargamma_n$ acts freely on $\mathfrak{h}$,
and the quotient space $ \bargamma_n \backslash \mathfrak{h}$ is a
punctured Riemann surface
with fundamental group~$\bargamma_n$.  If we supplement  $\mathfrak{h}$
with the projective rational line, $\mathbb{Q} \vee \{\infty\}$, then  we can think of the punctures as
cusps or $C_\infty$-points.  Then $\PSL_2(\integers_n)$
acts faithfully on the set of cusps  of $ \bargamma_n \backslash \mathfrak{h}$.

The following facts are well known.
\begin{enumerate}[\normalfont (1)]
\vskip-0.7cm \null
  \item  $\PSL_2(\integers_2)  = \gen{x, y \mid x^2 = y^3 = (xy)^2 = 1} = \Sym_3$, of order 6.
\vskip-0.7cm \null
  \item  $\barr( \bargamma_2) = 1$ and $ \bargamma_2$ is  free of rank two.
\vskip-0.7cm \null
  \item  $ \bargamma_2 \backslash \mathfrak{h}$ is a sphere with
three cusps, and $\PSL_2(\integers_2) \simeq \Sym_3$
acts naturally on the set of cusps.
\vskip-0.7cm \null
  \item  $\PSL_2(\integers_3)  = \gen{x, y \mid x^2 = y^3 = (xy)^3 = 1} = \Alt_4$, of order 12.
\vskip-0.7cm \null
  \item   $\barr( \bargamma_3) = 2$ and  $ \bargamma_3$ is  free of rank three.
\vskip-0.7cm \null
  \item
 $ \bargamma_3 \backslash \mathfrak{h}$ is a sphere with
four cusps, like a tetrahedron, and
$\PSL_2(\integers_3) \simeq \Alt_4$ acts naturally on the set of cusps.
\vskip-0.7cm \null
  \item  $\PSL_2(\integers_6)  \simeq \PSL_2(\integers_2) \times \PSL_2(\integers_3) \simeq \Sym_3 \times \Alt_4$, of order 72.
\vskip-0.7cm \null
  \item   $\barr( \bargamma_6) = 12$ and
$ \bargamma_6$ is  free of rank 13, and
$ \bargamma_6 = \bargamma_2 \cap \bargamma_3 $.
\vskip-0.7cm \null
  \item  $ \bargamma_6 \backslash \mathfrak{h}$ is a torus with
twelve cusps (see~\cite{CumminsPauli03}), and $\PSL_2(\integers_6) \simeq \Sym_3 \times \Alt_4$
acts faithfully on  the set of cusps.
\end{enumerate}
\vskip -0.8cm \hfill\qed \vskip 0.2cm
\end{example}

\begin{example}\label{ex:pp}  Let $p$ be an odd prime, and let $G = C_p \ast C_p$.

Then $\chi(G) = \frac{1}{p} + \frac{1}{p} - 1 = \frac{2-p}{p}$.

We have a presentation $G = \gen{x,y \mid x^p = y^p = 1}$.

In $\Sym_p$, consider  $x' = y' = x''  = (1,2,\ldots, p-1,p)$, and $y'' = \linebreak (p,p-1,\ldots,2,1)$.
The resulting kernels  $H$ and $K$ contain $x\overline y$ and $xy$, respectively;
recall that the overline indicates the inverse.
Now $HK$ contains $y \overline x$, $xy$ and $y^p$.
It follows that  $HK = G$.
By Proposition~\ref{prop:chars}, $\upp(G) \ge \frac{-1}{\chi(G)} = \frac{p}{p-2}$.
\hfill\qed
\end{example}

\begin{remark} Let us record  triples $(\barr(H), \barr(K), \barr(H\cap K))$
obtained in the above examples.
\begin{enumerate}[(i).]
\vskip-0.8cm \null
\item In $C_2\ast C_2 \ast C_2$, $(\barr(H), \barr(K), \barr(H\cap K)) = (1,2,4)$.
\vskip-0.8cm \null
\item In $C_2 \ast C_3$, $(\barr(H), \barr(K), \barr(H\cap K)) =  (1,2,12)$.
\vskip-0.8cm \null
\item In $C_2 \ast V$ and $C_2 \ast C_4$, $(\barr(H), \barr(K), \barr(H\cap K)) = (1,6,24)$.
\vskip-0.8cm \null
\item In $C_2 \ast C_p$, $p \ge 5$, $p$ prime, \newline$(\barr(H), \barr(K), \barr(H\cap K))
 = (p-2,\frac{1}{4}(p^2-1)(p-2),\frac{1}{4}(2p)(p^2-1)(p-2))$.
\vskip-0.8cm \null
\item In $C_p \ast C_p$, $p$ odd, $(\barr(H), \barr(K), \barr(H\cap K)) = (p-2, p-2, p(p-2))$. \hfill\qed
\end{enumerate}
\end{remark}

We now have a candidate for a sharp lower bound.

\begin{proposition}\label{prop:lower} If {\normalfont Notation~\ref{not:main}} holds,
then
$\upp(\mathcal{F}) \ge \depth(G){\cdot}\!\fun\!\height(G)$.
\end{proposition}

\begin{proof} Let $p = \height(G)$.

Thus
$p$ is $\infty$, $4$, or an odd prime, and $\fun\!\height(G) = \fun(p) = \frac{p}{p-2}$.

We consider two cases, with two subcases each.

\medskip

\noindent \textbf{Case 1.}   $\depth(G)=2$, that is, $G$ has an element of order two.

Here, there exists $j \in I$
such that $G_j$ has a subgroup which we can identify with~$C_2$.
By Remarks~\ref{rem:sec121}, we may assume that
$\abs{G_{\neg j}} \ge 3$.  Let $a$, $b$ and $c$ be three distinct elements of $G_{\neg j}$.

\textbf{Subcase 1.1.}  $p = \infty$.

We have
$C_2^a \ast C_2^b \ast C_2^c \,\,\le\,\,G_j^a \ast G_j^b \ast G_j^c \,\,\le\,\,G,$
and, hence,
$C_2 \ast C_2 \ast C_2$ embeds in $G$ in such a way that the finitely generated
free subgroups of $C_2 \ast C_2 \ast C_2$ are carried to  $\mathcal{F}$.

By Example~\ref{ex:222},
$\upp(\mathcal{F}) \ge 2 = 2{\cdot}\!\fun(\infty) = \depth(G){\cdot}\!\fun\!\height(G)$.

\textbf{Subcase 1.2.} $p$ is $4$ or an odd prime.

Here, there exists $i \in I$ such that $G_i$ has
a subgroup $P$ of order $p$.
Then $C_2^a \ast P^b \,\,\le\,\,G_j^a \ast G_i^b\,\,\le\,\,G,$
and, hence, $C_2 \ast P$ embeds in $G$ in such a way that the finitely generated
free subgroups of  $C_2 \ast P$ are carried to  $\mathcal{F}$.

By Examples~\ref{ex:2V} and~\ref{ex:2p},
$\upp(\mathcal{F}) \ge \frac{2p}{p-2} = 2{\cdot}\!\fun(p) = \depth(G){\cdot}\!\fun\!\height(G)$.

\medskip

\noindent\textbf{Case 2.} $\depth(G) = 1$, that is, $G$ has no element of order two.

\textbf{Subcase 2.1.}  $p  = \infty$.

In Notation~\ref{not:main}, we saw that
 $\upp(\mathcal{F}) \ge 1 = 1{\cdot}\!\fun(\infty) = \depth(G){\cdot}\!\fun\!\height(G)$.

\textbf{Subcase 2.2.}    $p$ is $4$ or an odd prime.

Notice that $p \ne 4$ since $\depth(G)\ne 2$.

Here,  there exists $j \in I$
such that $G_j$ has a subgroup which we can identify with~$C_p$.

By Remarks~\ref{rem:sec121}, $\abs{G_{\neg j}} \ge 2$.
Let $a$ and $b$ be two distinct elements of $G_{\neg j}$.
Then $C_p^a \ast C_p^b \,\,\le\,\, G_j^a \ast G_j^b \,\,\le\,\, G,$ and, hence,
$C_p \ast C_p$ embeds in $G$ in such a way that the finitely generated
free subgroups of $C_p \ast C_p$ are carried to  $\mathcal{F}$.

By Example~\ref{ex:pp},
$\upp(\mathcal{F}) \ge  \frac{p}{p - 2} = 1{\cdot}\fun(p) =  \depth(G){\cdot}\!\fun\!\height(G)$.
\end{proof}

\begin{exercise} Use the foregoing proof
to show that $\depth(G){\cdot}\!\fun\!\height(G)$~equals
\begin{equation}\label{eq:chars}
\max\{\textstyle\frac{(L:H)}{\barr(H)} \mid  H \in \mathcal{F},\,\, \barr(H) \ge 1, \,\, H \le L \le G, \,\,(L:H) < \infty \};
\end{equation}
here, $\frac{-1}{\chi(L)} = \frac{(L:H)}{\barr(H)}$.\hfill\qed
\end{exercise}

\section{Single-quotient subsets}\label{sec:3}

In this section, and in the next two sections,  $G$ is an arbitrary group.
Our main objective is to prove, in Corollary~\ref{cor:inequalities}(ii),
that, if  $A$ and $B$ are finite subsets of
 $G$ with at least two elements each, and  $\mathcal{C}$ is a set of pairwise-disjoint, single-quotient
subsets of $A \times B$, then
$\textstyle\sum\limits_{C \in \mathcal{C}} (\abs{C} - 2) \,\, \le \,\,
\fun\!\height(G){\cdot}(\abs{A} -2){\cdot}(\abs{B} -2).$   We recall that $\fun\!\height(G)$ was
described in Definitions~\ref{defs:begin}, and we now recall what we mean by a
`single-quotient' subset of $A\times B$.

\begin{definitions}\label{defs:var}   Let  $A$ and $B$ be finite subsets of $G$.

A subset $C$ of $A \times B$ is said to be a {\it single-product} subset of $A \times B$ if\linebreak
$\abs{\{ab \mid (a,b) \in C\}} = 1.$
Similarly,  $C$ is said to be a {\it single-quotient} subset if
$\abs{\{a \overline b \mid (a,b) \in C\}} = 1$.

For $x \in G$,
we let $\rep(x,A \times B)\,\, := \,\,\{(a,b) \in A \times B \mid ab = x\}\,\, \subseteq \,\, A \times B.$

For each positive integer $i$, we let
\begin{align*}
A{\cdot_i}B \quad&:=\quad \{x \in G : \abs{\rep(x,A \times B)} \ge i\} \quad \subseteq \quad G,\\
A{\cdot_{[=i]}}B \quad&:=\quad \{x \in G : \abs{\rep(x,A\times B)} = i\} \quad \subseteq \quad G.
\end{align*}
Thus, an element of $A{\cdot_i}B$, resp. $A{\cdot_{[=i]}}B$,
is an element of $G$ which has at least, resp. exactly, $i$ distinct representations
of the form $ab$ with $(a,b) \in A \times B$.

 We shall be interested in
$A{\cdot_1}B = AB$, $A{\cdot_2}B$, and $A{\cdot_{[=1]}}B = AB - A{\cdot_2}B$.
 \hfill\qed
\end{definitions}

The following result will be used frequently.

\begin{lemma}\label{lem:ii} For any finite subsets $A$, $B$ of $G$, the following hold.
\begin{enumerate}[\normalfont (i).]
\vskip-0.8cm \null
\item If $\abs{B} = 2$, then
$\abs{AB} + \abs{A{\cdot_2}B}\,\, = \,\,  2 \abs{A} +  2\abs{B} -4.$
\vskip-0.8cm \null
\item If $\abs{B}\ge 2$, then
$\abs{AB} + \abs{A{\cdot_2}B}\,\, \ge\,\,  2 \abs{A}.$
\end{enumerate}
\end{lemma}

\begin{proof}
Suppose that $b_1$ and $b_2$ are two distinct
elements of $B$, and let\linebreak $B' = \{b_1,b_2\}$.

 Then $B \supseteq B'$,
 $AB \supseteq AB' = Ab_1 \cup Ab_2$ and $A{\cdot_2}B \supseteq A{\cdot_2}B' = Ab_1 \cap Ab_2$.
Hence,
\begin{align*}
&\abs{AB} + \abs{A{\cdot_2}B} \ge \abs{AB'} + \abs{A{\cdot_2}B'} = \abs{Ab_1 \cup Ab_2} + \abs{Ab_1 \cap Ab_2}
= \abs{Ab_1} + \abs{Ab_2}\\&= 2\abs{A} = 2 \abs{A} +  2\abs{B'} -4.
\end{align*}
This proves (ii), and the case $B = B'$ proves (i).
\end{proof}

We call the next result {\it the key inequality}.
Recall from Definitions~\ref{defs:begin}
that $\height(G)$ is $\infty$ or $4$ or an odd prime, and that
$\fun\!\height(G) = \frac{\height(G)}{\height(G)-2} \in [1,3]$.

\begin{theorem}[= Theorem~\ref{th:presubmain}]\label{th:submain} For any finite subsets $A$, $B$ of $G$,
if $\abs{A} \ge 2$ and $\abs{B} \ge 2$, then
$\abs{AB} + \abs{A{\cdot_2}B} \,\, \ge \,\, \min\{ 2 \abs{A} + 2\abs{B} -4,\,\, 2{\cdot}\!\height(G)\}.$
\end{theorem}

\begin{proof} We postpone the lengthy proof to the next two sections; see Theorem~\ref{th:presubmain}.
\end{proof}

\begin{corollary}\label{cor:subsubmain} For any finite subsets $A$, $B$ of $G$, if $\abs{A} \ge 2$ and $\abs{B} \ge 2$, then
 $\abs{A} \abs{B} - \abs{AB} - \abs{A{\cdot_2}B} \,\,\,\, \le
 \,\,\,\, \fun\!\height(G){\cdot}(\abs{A}-2){\cdot}(\abs{B}-2).$
\end{corollary}

\begin{proof}   By symmetry, we may assume that $\abs{A} \ge \abs{B}$.

Let $p = \height(G)$.  Recall, from Definitions~\ref{defs:begin}, that
\begin{equation}\label{eq:decr}
\null\hskip -1.5cm \text{the function $\fun\colon [3,\infty] \to [1,3]$,  $x \mapsto \textstyle\frac{x}{x-2}$,
is strictly decreasing,}
\end{equation}
and $\fun\!\height(G) =\fun(p) = \frac{p}{p-2} \in [1,3]$.

We claim that at least one of the following holds.
\begin{enumerate}[(1).]
\vskip-0.7cm \null
\item $\abs{AB} + \abs{A{\cdot_2}B}\,\, \ge \,\, 2\abs{A} + 2\abs{B} - 4.$
\vskip-0.8cm \null
\item $\abs{A}\,\,  \ge\,\,  p.$
\vskip-0.8cm \null
\item $\abs{A} \,\, <\,\,  p
\quad\text{and}\quad \infty \,\, > \,\, \abs{AB} + \abs{A{\cdot_2}B} \,\, \ge\,\,  2p.$
\end{enumerate}
To see this, notice that if (1) and (2) fail, then (3) holds,
by Theorem~\ref{th:submain} (= Theorem~\ref{th:presubmain}).

 We now have three (overlapping) cases.

\medskip

\noindent\textbf{Case 1.}  $\abs{AB} + \abs{A{\cdot_2}B}\,\,  \ge \,\, 2\abs{A} + 2\abs{B} - 4$.
\vskip -1.2cm \null
\begin{align*}
\hskip -0.4cm \text{Here, }\abs{A} \abs{B} - \abs{AB} -\abs{A{\cdot_2}B}  \quad &\le\quad \abs{A} \abs{B} -2\abs{A} - 2\abs{B} + 4 \\
&=\quad \fun(\infty){\cdot}(\abs{A}-2){\cdot}(\abs{B}-2) \\
&\le \quad \fun(p){\cdot}(\abs{A}-2){\cdot}(\abs{B}-2) \text{ by~\eqref{eq:decr}.}
\end{align*}

\noindent\textbf{Case 2.}  $\abs{A}\,\,  \ge \,\, p\,\,  \ge \,\, 3$.
\vskip -1.2cm \null
\begin{align*}
\hskip -0.4cm \text{Here, }\abs{A} \abs{B} - \abs{AB} - \abs{A{\cdot_2}B}  \quad&\le
               \quad \abs{A} \abs{B} - 2\abs{A} \text{ by Lemma~\ref{lem:ii}(ii)}\\
& =\quad \textstyle \fun(\abs{A}){\cdot}(\abs{A}-2){\cdot}(\abs{B}-2) \\
&\le \quad \fun(p){\cdot}(\abs{A}-2){\cdot}(\abs{B}-2) \text{ by~\eqref{eq:decr}.}
\end{align*}

 \noindent\textbf{Case 3.}
$\abs{AB} + \abs{A{\cdot_2}B} \,\, \ge\,\,  2p\,\,$ and
$\infty \,\ > \,\, p \,\, >\,\,  \abs{A}\,\,  \ge \,\, \abs{B} \,\, \ge \,\, 2$.

Here,
\begin{align*}
(p -2)&(\abs{A}\abs{B} - \abs{AB} - \abs{A{\cdot_2}B})\quad\le \quad (p -2)(\abs{A}\abs{B} - 2p )  \\
&\le \quad (p -2)(\abs{A}\abs{B} - 2p ) + 2(p  - \abs{A})(p -\abs{B}) \\
&=\quad p\abs{A}\abs{B} - 2p^2 - 2\abs{A}\abs{B} + 4p + 2p^2 - 2p\abs{B} -2 p\abs{A} + 2\abs{A}\abs{B}\\
&=\quad p\abs{A}\abs{B}    + 4p  - 2p\abs{B} -2 p\abs{A} \quad \\&= \quad (p-2) {\cdot}\fun(p){\cdot}(\abs{A}-2){\cdot}(\abs{B}-2).
\end{align*}

\medskip

The desired result holds in all cases.
\end{proof}

Part (ii) of the following is the result that we shall apply in Section~\ref{sec:freeproduct}.

\begin{corollary}\label{cor:inequalities} Let  $A$ and $B$ be finite subsets of
a group $G$
such that $\abs{A} \ge 2$ and $\abs{B}\ge 2$.
\begin{enumerate}[\normalfont (i).]
\vskip-0.7cm \null
\item If $\mathcal{E}$ is a set of pairwise-disjoint, single-product
subsets of $A \times B$, then
$\textstyle\sum\limits_{E \in \mathcal{E}} (\abs{E} - 2) \,\, \le \,\,
\fun\!\height(G){\cdot}(\abs{A} -2){\cdot}(\abs{B} -2).$
\vskip-0.5cm \null
\item If $\mathcal{C}$ is a set of
pairwise-disjoint, single-quotient subsets of $A \times B$, then
$\textstyle\sum\limits_{C \in \mathcal{C}} (\abs{C} - 2) \,\, \le \,\,
\fun\!\height(G){\cdot}(\abs{A} -2){\cdot}(\abs{B} -2).$
\end{enumerate}
\end{corollary}

\begin{proof}   (i). If there exists some $E_0 \in \mathcal{E}$ such that
$\abs{E_0} \le 1$, then we may replace $\mathcal{E}$ with $\mathcal{E} -\{E_0\}$.
This respects the hypotheses and increases $\textstyle\sum\limits_{E \in \mathcal{E}} (\abs{E} - 2)$
by $2 - \abs{E_0}$.
By repeating this procedure as often as necessary, we may assume that, for each
$E \in \mathcal{E}$, $\abs{E} \ge 2$, and, hence, there exists a unique $x_E \in A{\cdot_2}B$ such that
 $\rep(x_E,A\times B) \supseteq E$.

If there exist some $E' \ne E''\in \mathcal{E}$ such that
$x_{E'} = x_{E''}$, then the disjoint union $E' \vee E''$ is again a single-product subset of $A \times B$,
and we may replace
$\mathcal{E}$ with $$\mathcal{E} -\{E',E''\} \cup \{E' \vee E''\}.$$
This respects the hypotheses and increases $\textstyle\sum\limits_{E \in \mathcal{E}} (\abs{E} - 2)$ by 2.
By repeating this procedure as often as necessary, we may assume that the map
$\mathcal{E} \to A{\cdot_2}B$, $E \mapsto x_E$, is injective.  Thus,
\begin{align*}
\textstyle\sum\limits_{E \in \mathcal{E}} (\abs{E} - 2)
 \,\, &\le \,\, \textstyle\sum\limits_{E \in \mathcal{E}}   (\abs{\rep(x_E,A\times B)}  - 2) \\
 \,\, &\le \,\, \textstyle \sum\limits_{x \in A{\cdot_2}B}  (\abs{\rep(x,A\times B)}  - 2)\\
&= \,\,\textstyle\abs{\bigvee\limits_{x \in A{\cdot_2}B}
\rep(x,A\times  B)} - 2 \abs{A{\cdot_2}B}\,\,\\&=\,\, \abs{A \times B}-\abs{A{\cdot_{[=1]}}B} - 2\abs{A{\cdot_2}B}\\
\,\,&= \,\, \abs{A} \abs{B}  - \abs{AB}- \abs{A{\cdot_2}B} \,\, \\&\le \,\,
\fun\!\height(G){\cdot}(\abs{A} -2){\cdot}(\abs{B} -2) \text{ by Corollary~\ref{cor:subsubmain}.}
\end{align*}

(ii). The bijection $A \times B \to A \times \overline B$, $(a,b) \mapsto (a,\overline b)$,
carries single-quotient subsets of $A \times B$ to single-product subsets of $A \times \overline B$.
Hence, by (i), we see that
$$\textstyle\sum\limits_{C \in \mathcal{C}} (\abs{C} - 2) \,\, \le \,\,
 \fun\!\height(G){\cdot}(\abs{A} -2){\cdot}(\abs{\overline B} -2),$$
and the result follows.
\end{proof}

\begin{examples} (i) Suppose that $G$ has an element $g$ whose order is at least~3.
Let $A = B = \{1,g,g^2\}$,
and let $\mathcal{C} = \{\{(1,1),(g, g),(g^2, g^2) \}\}$.
Here, Corollary~\ref{cor:inequalities}(ii) asserts that $(3-2) \le \fun\!\height(G){\cdot}(3-2)(3-2)$.

(ii). Suppose that $G$ has a finite, nontrivial subgroup $L$.
Let $A = B = L$, and let $\mathcal{C} = \{\{ (xy,y) \mid y \in L\} \mid x \in L\}$.  Here,
Corollary~\ref{cor:inequalities}(ii) asserts that
$\abs{L}{\cdot}(\abs{L} -2)  \,\, \le \,\, \fun\!\height(G){\cdot}(\abs{L} -2){\cdot}(\abs{L} -2)$.
\hfill\qed
\end{examples}

\section{Blocks and the Kemperman transform}\label{sec:4}

\begin{remarks}
To put the key inequality, Theorem~\ref{th:presubmain}/\ref{th:submain},
 into historical perspective, we record the following.\newline
\noindent \textbf{Kemperman's Theorem}. {\it
If  $A$ and $B$ are finite, nonempty subsets of a group~$G$, then
there exists a subgroup $L$ of $G$ such that
\begin{equation}\label{eq:kempe2}
\abs{A} + \abs{B} - \abs{AB}\,\,\le \,\, \abs{L} \,\,\le \,\,\abs{AB}.
\end{equation}
Moreover, if $A{\cdot_2}B \ne AB$, then  $\abs{L}$ can be taken
to be $1$}. \qed\newline This is a consequence of Theorems~5
and~3 of J.\ H.\ B.\ Kemperman's 1956 paper~\cite{Kemperman56};
it is a curious coincidence that 1956 also saw the publication
of Hanna Neumann's paper~\cite{HN1}. In the case where $G$ has
prime order,~\eqref{eq:kempe2} is the famous
Cauchy\d1Daven\-port Theorem, discovered by A.\
Cauchy~\cite{Cauchy13} in 1813 and by
{H.}~Davenport~\cite{Davenport35} in 1935.

We will be using (a variant of) the
marvellous
`Kemperman transform' which was introduced unnamed in the
proofs of Theorems~5 and~3 of~\cite{Kemperman56}; see
Definition~\ref{def:bde}, below. Kemperman pointed out that
this transform is closely related to the type of reasoning that
H.\ B.\ Mann~\cite{Mann} had employed  to prove the
Lan\-dau--Schur--Khint\-chine $\alpha+\beta$-conjecture.
\hfill\qed
\end{remarks}

\medskip

In this section, we introduce concepts that will be used in the
proof in the next section.

\begin{definitions}\label{defs:sound}
For each $n \in \naturals$, we let
 $\mathcal{S}_n$ denote the set of pairs $(A,B)$ such that $A$ and $B$ are
finite subsets of $G$ with $\abs{A} \ge n$ and $\abs{B} \ge n$.  We shall be interested in
 $\mathcal{S}_2 \subseteq \mathcal{S}_0$.

For $(A,B) \in \mathcal{S}_0$, we define
$\hbar(A,B) := \abs{AB} + \abs{A{\cdot_2}B} - 2\abs{A} - 2\abs{B} \in \integers$.

By a {\it block} (in $G$) we mean a subset of $G$ of the form
$cPd$ where $c$ and $d$ are elements of $G$, and $P$ is a
subgroup of $G$ whose order is either $4$ or an odd prime. We
remark that $\abs{cPd} = \abs{P} \ge \height(G)$. By replacing
the triple $(c,P,d)$ with the triple $(cd,P^d,1)$, we can
arrange that $d=1$.

If $C$ is a finite subset of $G$, we let $\blocks(C)$ denote the number of subsets
of $C$ which are blocks in $G$.

An element $(A,B)$  of $\mathcal{S}_2$ is said to be {\it sound} if
(at least) one of the following holds: $\abs{AB} + \abs{A{\cdot_2}B}\,\,\ge\,\,2\abs{A} + 2\abs{B} -4$
(equivalently, $\hbar(A,B) \ge -4$),
or $\blocks(A{\cdot_2}B) \ge 1$, or $\blocks(AB) \ge 2$.
\hfill\qed
\end{definitions}

In the next section, we shall show that every element of $\mathcal{S}_2$ is sound.

\begin{examples} (i).  Suppose that $G$ has an element $g$ whose order is at least 3, and
take $A = B = \{1,g\}$.

 Then $AB = \{1,g,g^2\}$ and $A{\cdot_2}B = \{g\}$.

Here, $\hbar(A,B) = \abs{AB} + \abs{A{\cdot_2}B} - 2\abs{A} - 2\abs{B}
= 3 + 1 -4 - 4 = -4$.

Also,
$\blocks(A{\cdot_2}B) = 0$ and $\blocks(AB) \le 1$.

(ii). Suppose that $G$ has a subgroup $P$ of order $4$ or an odd prime, and take $A = B = P$.

Then $AB = A{\cdot_2}B = P$.

Here, $\blocks(AB) = \blocks(A{\cdot_2}B) = 1$.

Also, $\hbar(A,B) = \abs{AB} + \abs{A{\cdot_2}B} - 2\abs{A} - 2\abs{B}
= - 2\abs{P} < -4.$

 (iii).  We do not know of an example where $\blocks(AB) \ge 2$ but\linebreak
$\hbar(A,B) < -4$ and $\blocks(A{\cdot_2}B) = 0$. \hfill\qed
\end{examples}

\begin{lemma}\label{lem:sound}  Let $(A,B) \in \mathcal{S}_2$.  If $(A,B)$ is sound, then
\begin{equation}\label{eq:main}
\abs{AB} + \abs{A{\cdot_2}B} \quad \ge \quad \min\{ 2 \abs{A} + 2\abs{B} -4,\,\, 2{\cdot}\!\height(G)\}.
\end{equation}
\end{lemma}

\begin{proof}
{From} Definitions~\ref{defs:sound}, we have three possibilities.

\medskip

\noindent\textbf{Case 1.}  $\abs{AB} + \abs{A{\cdot_2}B}\,\, \,\ge  \,\, \, 2\abs{A} + 2\abs{B} -4$.

Here,~\eqref{eq:main} holds.

\medskip

\noindent\textbf{Case 2.}  $\blocks(A{\cdot_2}B) \ge 1$.

Here, $\abs{A{\cdot_2}B} \ge \height(G)$.  Hence,
$\abs{AB} + \abs{A{\cdot_2}B}\,\,\ge\,\, 2{\cdot}\abs{A{\cdot_2}B}\,\, \ge\,\, 2{\cdot}\!\height(G),$
and, hence,~\eqref{eq:main} holds.

\medskip

\noindent\textbf{Case 3.}  $\blocks(AB) \ge 2$.

We subdivide this case into two subcases.

\textbf{Subcase 3.1.} \,\, $\abs{A{\cdot_2}B} \ge 2$.

We have $AB \supseteq c_1P_1 \cup c_2P_2$
where $c_1P_1$ and $c_2P_2$ are two different blocks in~$G$.

We claim that $\abs{c_1P_1 \cap c_2P_2} \le 2$.
  Suppose that
$d$ is an element of\linebreak $c_1P_1 \cap c_2P_2$.  Then
$dP_1 = c_1P_1$ and $dP_2 = c_2P_2$. Hence, $dP_1 \ne dP_2$, and, hence,
$P_1 \ne P_2$, and, hence,
 $\abs{P_1 \cap P_2} \le 2$, by the conditions on the orders.  Now,
 $c_1P_1 \cap c_2P_2 =  dP_1 \cap dP_2 = d(P_1\cap P_2),$
and the claim is proved.

Thus $\abs{AB} \ge \abs{c_1P_1} + \abs{c_2P_2} - \abs{c_1P_1 \cap c_2P_2}
\ge \height(G) + \height(G) -2.$

Since $\abs{A{\cdot_2}B} \ge 2$, we see that
$\abs{AB} + \abs{A{\cdot_2}B} \ge 2{\cdot}\!\height(G)$, and~\eqref{eq:main} holds.

\textbf{Subcase 3.2.} \,\, $\abs{A{\cdot_2}B} \le 1$.

 If $\abs{B} = 2$, then~\eqref{eq:main} holds by Lemma~\ref{lem:ii}(i).
Thus, we may assume that
$\abs{B} \ge 3$.
Here,
$(\abs{A} - 2){\cdot}(\abs{B}-2 - \abs{A{\cdot_2}B)} \ge (\abs{A} - 2){\cdot}(3-2-1) = 0$, and it follows that
\begin{equation}\label{eq:adhoc}
\abs{A} {\cdot}(\abs{B} - \abs{A{\cdot_2}B}) + 2\abs{A{\cdot_2}B}  \ge 2\abs{A} + 2\abs{B} -4.
\end{equation}

Since $A{\cdot_{[=1]}}B
= \textstyle\bigvee\limits_{a \in A} (aB \cap A{\cdot_{[=1]}}B)  =
\bigvee\limits_{a \in A} (aB - A{\cdot_2}B)$, we see that
\begin{align*}
\\[-0.8cm]
&\abs{A{\cdot_{[=1]}}B} = \abs{\textstyle\bigvee\limits_{a \in A} (aB - A{\cdot_2}B)} =
\textstyle\sum\limits_{a \in A} \abs{aB - A{\cdot_2}B} \\&\ge
\textstyle\sum\limits_{a \in A} (\abs{aB} - \abs{A{\cdot_2}B}) =
\textstyle\sum\limits_{a \in A} (\abs{B} - \abs{A{\cdot_2}B}) = \abs{A} {\cdot}(\abs{B} - \abs{A{\cdot_2}B}).
\intertext{Now,}\\[-1.1cm]
&\abs{AB} + \abs{A{\cdot_2}B} = \abs{A{\cdot_{[=1]}}B} + 2\abs{A{\cdot_2}B} \\
&\ge
\abs{A} {\cdot}(\abs{B} - \abs{A{\cdot_2}B}) + 2\abs{A{\cdot_2}B} \text{ by the foregoing}
\\&\ge 2\abs{A} + 2\abs{B} -4 \text{ by~\eqref{eq:adhoc}},
 \end{align*}
  and, hence,~\eqref{eq:main} holds.

\medskip

Thus,~\eqref{eq:main} holds in all cases.
\end{proof}

\begin{definitions}   We endow $\mathcal{S}_0$ with a partial order by
assigning  four indicators to each $(A,B) \in \mathcal{S}_0$.

 The {\it first indicator} of $(A,B)$ is $\abs{AB} \in \naturals \subset \integers$.

The {\it second indicator}  of $(A,B)$ is $\hbar(A,B) = \abs{AB} + \abs{A{\cdot_2}B} - 2\abs{A} - 2\abs{B} \in \integers.$

The {\it third indicator} of $(A,B)$ is  $\abs{B} \in \naturals \subset \integers$.

The {\it fourth indicator} of $(A,B)$ is $\abs{A} \in \naturals \subset \integers$.

We say that the {\it indicator sequence} of $(A,B)$ is $(\abs{AB}, \hbar(A,B), \abs{B}, \abs{A})$.

Considered lexicographically, the indicator sequence gives a partial order,
 denoted $\succcurlyeq$, on~$\mathcal{S}_0$.
Thus, if  $(A',B')$ is an element of $\mathcal{S}_0$, we write $(A,B) \succ  (A',B')$ if
and only if
$$( \abs{AB},  \hbar(A,B), \abs{B}, \abs{A})\,\, > \,\, ( \abs{A'B'},\hbar(A',B'), \abs{B'}, \abs{A'})$$
in the lexicographic ordering of $\integers^4$.
\hfill \qed
\end{definitions}

\begin{lemma} \label{lem:well} There are
no infinite, strictly descending chains in $(\mathcal{S}_2, \succcurlyeq)$.
\end{lemma}

\begin{proof} Recall that the indicator sequence of $(A,B)$ is $(\abs{AB}, \hbar(A,B), \abs{B}, \abs{A})$.
 In any infinite descending chain in $(\mathcal{S}_2, \succcurlyeq)$,
the first  indicator eventually becomes constant.
Once the first indicator is constant, the other three indicators can take only finitely many values, and,
hence, eventually become constant also.

This is also true in $\mathcal{S}_1$,
but not in $\mathcal{S}_0$.
\end{proof}

\begin{notation}\label{not:delta} Let us think of
$\{A,B,\cdot_1, \cdot_2, \hbar\}$ as a set of five functions with domain~$\mathcal{S}_0$,
where  $A$ and $B$ denote the projections onto the
first and second coordinates, respectively, of elements of  $\mathcal{S}_0$.

Let $(A_1,B_1),\,(A_2,B_2) $ be elements of $ \mathcal{S}_0$.

 We define a map $\delta = \delta^{((A_2,B_2),(A_1,B_1))} \colon\{A,B,\cdot_1, \cdot_2, \hbar\} \to \integers$
with the following values:
\vskip-0.7cm \null
\begin{list}{}
\item $\delta(A):= \abs{A_2} - \abs{A_1}$; \quad  $\delta(B):= \abs{B_2} -\abs{B_1}$;
\item $\delta ( {\cdot_1}) := \abs{A_2B_2}- \abs{A_1B_1}$;\quad
$\delta ( {\cdot_2}) := \abs{(A_2){\cdot_2}(B_2)}- \abs{(A_1){\cdot_2}(B_1)}$;
\vskip-0.7cm \null
\item $\delta (\hbar):=  \hbar(A_2,B_2) - \hbar(A_1,B_1) = \delta({\cdot_1}) +
\delta({\cdot_2}) - 2\delta(A) - 2\delta(B)$. \vskip 0.3cm
\end{list}
 In applications, $A_1$ will always be denoted $A$, with little risk of
confusion.  \hfill\qed
\end{notation}

\begin{definition}\label{def:bde}  Let $(A,B) \in \mathcal{S}_0$ and let
$x \in G$.

Set  $(A^{+},B^{-}) = (A \cup Ax, B \cap \overline{x}B)$ and
$(A^{-},B^{+}) = (A \cap A \overline x, B \cup {x}B)$.
Clearly,
\begin{equation}\label{eq:contains}
A^{+}B^{-} \subseteq AB \text{\,\,\, and \,\,\,} A^{-}B^{+} \subseteq AB.
\end{equation}

With Notation~\ref{not:delta}, let $\delta^{+} = \delta^{((A^{+},B^{-}), (A,B))}$ and
$\delta^{-} = \delta^{((A^{-},B^{+}), (A,B))}$.

We define the (revised) {\it Kemperman transform of $(A,B)$ with respect to $x$} to be
\begin{align}\label{eq:kempe}
(A',B'):&=\begin{cases}
 (A^{-},B^{+})  &\text{if $\delta^{-}(\hbar) < 0$,}\\
 (A^{+},B^{-})&\text{if $\delta^{-}(\hbar) \ge 0$ and $\delta^{+}(\hbar) < 0$,}\\
 (A^{+},B^{-})&\text{if $\delta^{-}(\hbar)  \ge 0$ and $\delta^{+}(\hbar)\ge 0$
and $\delta^{+}(B) < 0$,}\\
 (A^{-},B^{+}) &\text{if $\delta^{-}(\hbar)  \ge 0$ and $\delta^{+}(\hbar)\ge 0$
and $\delta^{+}(B) \ge 0$.}\\
\end{cases}
\end{align}
Thus $(A',B')$ is a well-defined element of $\mathcal{S}_0$.

We now make a sequence of remarks about this construction.

We call the bijection $G \times G \to G \times G$, $(a,b) \mapsto (\overline b, \overline a)$,
the {\it dual} map.
  Any statement about $G \times G$ can be ``dualized'' in a natural way.

\Lsub{}\label{comm:1} \hskip -0.25cm \textbf{Remark.}
\quad $\delta^{+}(A) + \delta^{-}(A) = \delta^{+}(B) + \delta^{-}(B) = 0.$

\begin{proof}
Notice that $\abs{A - A\overline x} = \abs{(Ax - A) \overline x)} =  \abs{Ax - A}$.  Now,
\begin{align*}\label{eq:deltaA}
&\delta^{+}(A)+\delta^{-}(A) \,\,= \,\,(\abs{A^+} - \abs{A}) + (\abs{A^-} - \abs{A})\\ \nonumber
&=\,\,(\abs{A \cup A x} - \abs{A}) + (\abs{A \cap A\overline x} - \abs{A})
\,\, = \,\,\abs{A x - A} - \abs{A - A\overline x} = 0.
\end{align*}
Dualizing, we see that $\delta^{+}(B) + \delta^{-}(B) = 0.$
\end{proof}

\Lsub{}\label{comm:2} \hskip -0.25cm \textbf{Remark.}
\quad
$\delta^{+}(\cdot_1)\,\, =\,\, -\abs{AB - A^+ B^-} \,\,\, \le \,\, 0$ and\newline
$\delta^{-}(\cdot_1)\,\, =\,\, -\abs{AB - A^- B^+} \,\, \, \le \,\, 0$.

\begin{proof}  This is clear from~\eqref{eq:contains}.
\end{proof}

\Lsub{}\label{comm:3} \hskip -0.25cm \textbf{Remark.}
\quad
$\max\{0, \delta^-(\cdot_2)\} \quad \le \quad \abs{(A^-){\cdot_2}(B^+)\,\,\cap\,\,A{\cdot_{[=1]}}B}$.

\begin{proof}
\begin{align*}
\hskip 2.3cm  \\[-1.55cm]
\delta^{-}(\cdot_2) \,\,&=\,\, \abs{(A^-){\cdot_2}(B^+)} - \abs{A{\cdot_2}B}  \\ \nonumber
&= \,\, \abs{(A^-){\cdot_2}(B^+)\,\,\cap\,\,AB}    -  \abs{A{\cdot_2}B} \text{ since $A^-B^+ \subseteq AB$}\\ \nonumber
&= \,\, \abs{(A^-){\cdot_2}(B^+)\,\,\cap\,\,A{\cdot_{[=1]}}B} +
 \abs{(A^-){\cdot_2}(B^+)\,\,\cap\,\,A{\cdot_{2}B}}  -  \abs{A{\cdot_2}B} \\ \nonumber
&\le\,\,\abs{(A^-){\cdot_2}(B^+)\,\,\cap\,\,A{\cdot_{[=1]}}B}.\qedhere
\end{align*}
\end{proof}

\Lsub{}\label{comm:4} \hskip -0.25cm \textbf{Remark.}
\quad $A^+ B^- \,\,\cap\,\,(A^-){\cdot_2}(B^+) \,\, \cap\,\,\ A{\cdot_{[=1]}}B \quad = \quad \emptyset.$

\begin{proof} Suppose that
\begin{equation}\label{eq:imposs}
c \in  A^+ B^- \,\,\cap\,\,(A^-){\cdot_2}(B^+) \,\, \cap\,\,\ A{\cdot_{[=1]}}B,
\end{equation}
and let $(a,b)$ denote the unique element of $\rep(c,A \times B)$.

By~\eqref{eq:imposs}, the equation $c = a'b'$  has at least two solutions
$(a',b')$ with $(a',b')$ in $(A^-)\times(B^+) = (A \cap A\overline x) \times (B \cup  x B)$.

\begin{enumerate}[\text{Type} 1.]
\vskip-0.6cm \null
  \item
$b' \in B$. \newline
Here, $(a',b') \in \rep(c,A \times B) = \{(a,b)\}$.   Hence, $a = a' \in A \cap A\overline x$.
Observe that if $b \in xB$ then
$(ax, \overline xb) \in \rep(ab,A \times B) = \{(a,b)\}$, which is a contradiction;
hence, here, $b \in B - xB$.
\vskip-0.6cm \null
  \item $b' \in xB - B$.\newline
Here, $(a'x,\overline  xb') \in  \rep(c,A \times B) = \{(a,b)\}$.
Hence, $a = a'x \in Ax$.
Moreover,  $b = \overline xb' \in B - \overline xB = B - B^-$.
Here, $(a',b') = (a\overline  x,xb)$.
\vskip-0.6cm \null
\end{enumerate}
In summary, the equation $c = a'b'$ has {\it exactly} two solutions
 $(a',b') \in A^- \times B^+$, one of each type, namely, $(a,b)$ and $(a\overline  x, x b)$.

It follows that $(a,b) \in (Ax \cap A \cap A\overline x) \times (B - (\overline x B \cup xB))$.

By~\eqref{eq:imposs}, there exists some $(a'',b'') \in A^+ \times B^- = (A \cup Ax) \times (B \cap \overline xB)$
such that  $a''b'' = c$.

\begin{enumerate}[\text{Case }1.]
\vskip-0.8cm \null
  \item $a'' \in A$. \newline
Here, $(a'',b'')\in  \rep(c,A \times B) = \{(a,b)\}$. Hence, $b = b'' \in B^- \subseteq \overline xB$.
This contradicts the fact that $b \in B - \overline xB$.
\vskip-0.8cm \null
  \item $a'' \in Ax - A$. \newline
Here, $(a''\overline x, x b'') \in\rep(c,A \times B) = \{(a,b)\}$.
Hence, $a = a''\overline x   \in A - A\overline x$.
This contradicts the fact that $a \in A\overline x$.
\end{enumerate}
\vskip-0.6cm \null
This completes the proof of Remark~\ref{comm:4}.
\end{proof}

On dualizing Remark~\ref{comm:4}, we get the following.

\Lsub{}\label{comm:5} \hskip -0.25cm \textbf{Remark.}
\quad $A^- B^+ \,\, \cap \,\, (A^+){\cdot_2}(B^-) \,\, \cap\,\,\ A{\cdot_{[=1]}}B\quad = \quad \emptyset.$
\hfill\qed

\Lsub{}\label{comm:6} \hskip -0.25cm \textbf{Remark.}
\quad $\delta^{+}(\hbar) + \delta^{-}(\hbar) \le 0.$

\begin{proof}
 Here,
\begin{align}\label{eq:deltaminus}
\delta^{-}(\cdot_2) \,\, &\le\,\,\abs{(A^-){\cdot_2}(B^+)\,\,\cap\,\,A{\cdot_{[=1]}}B}  \text{ by Remark~\ref{comm:3}}\\ \nonumber
\,\,&\le\,\, \abs{AB - A^+ B^-} \text{ by Remark~\ref{comm:4}} \\ \nonumber
&=\,\, -\delta^{+}(\cdot_1) \text{ by Remark~\ref{comm:2}}.
\end{align}
By dualizing, we see that
\begin{equation}\label{eq:otherdelta}
\delta^{+}(\cdot_2) \,\,\le\,\,-\delta^{-}(\cdot_1).
\end{equation}

By combining Remark~\ref{comm:1} with~\eqref{eq:deltaminus} and~\eqref{eq:otherdelta}, we obtain
\begin{align*}
\delta^{+}(\cdot_1) + \delta^{+}(\cdot_2) - &2 \delta^{+}(A)
- 2 \delta^{+}(B)\\ \nonumber
& \qquad +\delta^{-}(\cdot_1) + \delta^{-}(\cdot_2) - 2 \delta^{-}(A)
- 2 \delta^{-}(B)) \quad \le \quad 0,
\end{align*}
and Remark~\ref{comm:6} is proved.
\end{proof}

\Lsub{}\label{comm:7} \hskip -0.25cm \textbf{Remark.}  The following holds:
\begin{equation*}
(A',B') =\begin{cases}
 (A^{-},B^{+})  &\text{if $\delta^{-}(\hbar) < 0$,}\\
 (A^{+},B^{-})&\text{if $\delta^{-}(\hbar) \ge 0$ and $\delta^{+}(\hbar) < 0$,}\\
 (A^{+},B^{-})&\text{if $\delta^{-}(\hbar)  = \delta^{+}(\hbar)= 0$
and $\delta^{+}(B) < 0$,}\\
 (A^{-},B^{+}) &\text{if $\delta^{-}(\hbar)  = \delta^{+}(\hbar)= 0$
and $\delta^{+}(B) = 0$.}\\
\end{cases}
\end{equation*}
Of course, if $(A^{+},B^{-})  = (A^{-},B^{+})$, then
$(A^{+},B^{-})  = (A^{-},B^{+}) = (A,B)$.

\begin{proof}
The description of  $(A',B')$ follows from~\eqref{eq:kempe}, and
Remark~\ref{comm:6}, and the fact that  $\delta^{+}(B)\le 0$; recall that
$\delta^{+}(B) = \abs{B^-} - \abs{B}$.
\end{proof}

This completes the desired description of the Kemperman transform.\hfill\qed
\end{definition}

\section{Proof of the key inequality}\label{sec:5}

This section is structured as the proof of the key inequality.
Recall Definitions~\ref{defs:sound}.  We fix, throughout the
proof, an element $(A,B)$ of~$\mathcal{S}_2$ and we show that
$(A,B)$ is sound by progressively finding various assumptions
that we are free to make.

\begin{assumptions}\label{ass:key}  Let $(A,B)$ be an element of~$\mathcal{S}_2$.
We want to show that $(A,B)$ is sound.

By Lemma~\ref{lem:well} and transfinite induction,
we have the following (transfinite) induction hypothesis:
we assume, without loss of generality, that
 in $(\mathcal{S}_2, \succcurlyeq)$,
 every element  which
is strictly smaller than $(A,B)$ is sound.
\hfill\qed
\end{assumptions}

\begin{lemma}  With {\normalfont Assumptions~\ref{ass:key}}, if $\abs{A} < \abs{B}$, then $(A,B)$ is sound.
\end{lemma}

\begin{proof} Recall that the indicator sequence of $(A,B)$ is $(\abs{AB}, \hbar(A,B), \abs{B}, \abs{A})$.
In passing from $(A,B)$ to its dual, $(\overline B, \overline A)$,  the first two indicators stay the same,
while the third indicator decreases by $\abs{B} -\abs{A}$.  By the induction hypothesis,
 Assumptions~\ref{ass:key},   $(\overline B, \overline A)$ is sound.  Dualizing, we see that
$(A,B)$ is sound.
\end{proof}

Also, by Lemma~\ref{lem:ii}(i),  $(A,B)$ is sound if $\abs{B} = 2$.

\begin{assumptions}\label{ass:ge} We assume,
without loss of generality, that $\abs{A} \ge \abs{B} \ge 3$.  \hfill\qed
\end{assumptions}

\begin{lemma}\label{lem:cols}   With {\normalfont Assumptions~\ref{ass:key} {\it and}~\ref{ass:ge}}, the following hold.
\begin{enumerate}[\normalfont (i).]
\vskip -0.6cm\null
\item If, for some $a \in A$, $\abs{aB\,\, \cap\,\, A{\cdot_{[=1]}}B} \ge 2$, then $(A,B)$ is sound.
\vskip -0.6cm\null
\item If,  for some $b \in B$,   $\abs{Ab\,\, \cap\,\, A{\cdot_{[=1]}}B} \ge 2$, then $(A,B)$ is sound.
\end{enumerate}
\end{lemma}

\begin{proof}
For (i), set $(A',B') = (A-\{a\},B)$; for (ii), set $(A',B') = (A,B-\{b\})$.

In both cases, $(A',B') \in \mathcal{S}_2$, by Assumptions~\ref{ass:ge}.

It is easy to see that, for (i), $A'B' =   AB - (aB \cap A{\cdot_{[=1]}}B)$, while, for (ii),
$A'B' =   AB - (Ab \cap A{\cdot_{[=1]}}B)$.

In both cases, $A'{\cdot_2}B'  \subseteq  A{\cdot_2}B$.

Thus, in both cases, $\abs{A'} + \abs{B'} = \abs{A} + \abs{B} -1$,
$\abs{A'B'} \le \abs{AB} -2$, and $\abs{A'{\cdot_2}B'} \le \abs{A{\cdot_2}B}$.
Now the two cases are handled together.

Recall that the indicator sequence of $(A,B)$ is $(\abs{AB}, \hbar(A,B), \abs{B}, \abs{A})$.
In passing from $(A,B)$ to $(A',B')$,  the first indicator decreases, by at least~2.
By the induction hypothesis, Assumptions~\ref{ass:key}, $(A',B')$ is sound.
By Definitions~\ref{defs:sound}, there are three possibilities.

\medskip

\noindent\textbf{Case 1.} $\abs{A'B'} + \abs{A'{\cdot_2}B'}\,\, \,\ge  \,\, \, 2\abs{A'} + 2\abs{B'} -4$.

Here,\newline
$
\abs{AB} + \abs{A{\cdot_2}B} \,\,
\ge \,\, 2 + \abs{A'B'} + \abs{A'{\cdot_2}B'} \,\,  \ge \,\, 2 +  2\abs{A'} + 2\abs{B'} -4 \,\,  = \,\,
 2\abs{A} + 2\abs{B} -4.
$
Thus, $(A,B)$ is sound.

\medskip

\noindent\textbf{Case 2.}  $\blocks(A'{\cdot_2}B') \ge 1$.

Since $A{\cdot_2}B \supseteq A'{\cdot_2}B'$, we see that $\blocks(A{\cdot_2}B) \ge \blocks(A'{\cdot_2}B') \ge 1$,
and $(A,B)$ is sound.

\medskip

\noindent\textbf{Case 3.}   $\blocks(A'B') \ge 2$.

Since $AB \supseteq A'B'$,  we see that $\blocks(AB) \ge \blocks(A'B') \ge 2$, and $(A,B)$ is sound.

\medskip

Hence, (i) and (ii) hold.
\end{proof}

\begin{assumptions}\label{ass:Q1cols} We assume, without loss of generality, that the following hold.
\begin{enumerate}[\rm (i).]
\vskip-0.7cm \null
  \item  For each $a \in A$, $\abs{aB \,\, \cap\,\, A{\cdot_{[=1]}}B}\,\, \le \,\,1$.
\vskip-0.7cm \null
\item   For each $b \in B$, $\abs{Ab \,\, \cap\,\, A{\cdot_{[=1]}}B}\,\, \le \,\,1$.
\end{enumerate}
\vskip -0.8cm \hfill\qed \vskip 0.2cm
\end{assumptions}

The proofs of Lemmas~\ref{lem:e} and~\ref{lem:bd}, which are
modelled on the proofs of Theorem~5 and Theorem~3
of~\cite{Kemperman56}, respectively, have a large common part
which we now describe.

\begin{hypotheses}\label{hyp} With Assumptions~\ref{ass:key}, ~\ref{ass:ge}, and~\ref{ass:Q1cols},
let $x$ be an element of $G$ such that $A \ne Ax$ and let
 $(A',B')$ be the  Kemperman transform
of $(A,B)$ with respect to $x$, with
notation as in Definition~\ref{def:bde}.

Since $Ax \ne A$, we see that $x \ne 1$, and that $A^- \subset A \subset A^+$.

\Lsub{}\label{cons:1} \hskip -0.25cm \textbf{Consequence.}
\quad {\it If \,$1, x \in B$, and
$$2\abs{A^+} + \abs{(A^-){\cdot_2}(B^+)\,\,\cap\,\,A{\cdot_{[=1]}}B}\quad
\ge\quad 2\abs{A}+2\abs{B} -2,$$
then $(A,B)$ is  sound.}

\begin{proof}
Observe that $AB \cap A^+B^- \,\,\supseteq \,\,(A1 \cup Ax) \cap (A^+1)\,\, =\,\, A^+$.
Hence, by Remark~\ref{comm:4},
$$\abs{AB} \quad \ge \quad \abs{A^+} + \abs{(A^-){\cdot_2}(B^+)\,\,\cap\,\,A{\cdot_{[=1]}}B}.$$
Since $A{\cdot_2}B \,\, \supseteq \,\, A^+ -  (A{\cdot_{[=1]}}B),$
it follows from Assumptions~\ref{ass:Q1cols}(ii) that $$\abs{A{\cdot_2}B}  \ge   \abs{A^+} -2.$$
Hence, $$\abs{AB} + \abs{A{\cdot_2}B}\, \ge
 \, 2\abs{A^+} + \abs{(A^-){\cdot_2}(B^+)\,\,\cap\,\,A{\cdot_{[=1]}}B} -2\,\, \,\ge  \,\, \, 2\abs{A} + 2\abs{B} -4,$$
and $(A,B)$ is sound.
\end{proof}

\Lsub{}\label{cons:2} \hskip -0.25cm \textbf{Consequence.}
\quad {\it If
$$ 2\abs{A^+} + \abs{(A^-){\cdot_2}(B^+)\,\,\cap\,\,A{\cdot_{[=1]}}B}\quad \le\quad 2\abs{A}+2\abs{B} -3,$$
then $(A',B') \in \mathcal{S}_2$.}

\begin{proof} By Remark~\ref{comm:3}, the hypothesis implies that
\begin{align}
2\abs{A^+} + 0\quad &\le\quad 2\abs{A}+2\abs{B} -3, \text{ and,}\label{eq:ineq1}\\
2\abs{A^+} + \delta^-(\cdot_2)\quad &\le\quad 2\abs{A}+2\abs{B} -3.\label{eq:ineq2}
\end{align}
\begin{enumerate}[\text{Case} 1.]
\vskip-0.8cm \null
  \item $(A',B') = (A^-, B^+)$.

Using~\eqref{eq:ineq1} and Assumptions~\ref{ass:ge}, we see that
$$2 \abs{A^-} \,\,= \,\,2(2\abs{A} - \abs{A^+})
\,\,\ge \,\,2\abs{A} - 2\abs{B} + 3  \,\,\ge\,\, 0 + 3.$$
Thus, $\abs{A^-} \ge \frac{3}{2}$, and, hence, $(A',B')  \in \mathcal{S}_2$.

\vskip-0.7cm \null
  \item $(A',B') = (A^+,B^-)$.

It follows from Remark~\ref{comm:7} that $\delta^{-}(\hbar) \ge
0$.  Hence
\begin{align*}
0 \quad &\le \quad \delta^{-}(\hbar) \quad\\&=\quad \delta^{-}(\cdot_1) + \delta^{-}(\cdot_2)
- 2 \delta^{-}(A) - 2\delta^{-}(B)\\
&\le\quad  0 + \delta^{-}(\cdot_2)
- 2 \delta^{-}(A) - 2\delta^{-}(B) \text{\,\,\,\,\,by Remark~\ref{comm:2}}\\
\quad&=\quad  \delta^{-}(\cdot_2) + 2 \delta^{+}(A)+ 2\delta^{+}(B)
                 \text{\,\,\,\,\,by Remark~\ref{comm:1}}\\
&=\quad \delta^{-}(\cdot_2) + 2\abs{A^+} - 2\abs{A} + 2\abs{B^-} - 2\abs{B}\\
\quad&\le\quad - 3 + 2 \abs{B^-} \text{\quad by~\eqref{eq:ineq2}}.
\end{align*}
Here, $\abs{B^-} \ge \frac{3}{2}$, and, hence, $(A',B') \in \mathcal{S}_2$.
\end{enumerate}

\noindent In all cases then, $(A',B')  \in \mathcal{S}_2$ and Consequence~\ref{cons:2} is proved.
 \end{proof}

\Lsub{}\label{conss:3} \hskip -0.25cm \textbf{Consequences.}
 {\it The following hold: $A'B' \subseteq AB$;
 $\hbar(A',B') \le \hbar(A,B)$;    $(A,B) \succ (A',B')$;
and, if $(A',B') \in \mathcal{S}_2$, then $(A',B')$ is sound.
}

\begin{proof} The first assertion follows from~\eqref{eq:contains}.

With Notation~\ref{not:delta}, let $\delta' = \delta^{((A',B'), (A,B))}$.
  It follows from Remark~\ref{comm:7} that
$\delta'(\hbar) \le 0$, and, hence, $\hbar(A',B') \le
\hbar(A,B)$.

Recall that the indicator sequence of $(A,B)$ is $(\abs{AB}, \hbar(A,B), \abs{B}, \abs{A})$.
 We now discuss how the four indicators
change in passing from $(A,B)$ to $(A',B')$.  We have just seen that
the first two indicators do not increase.

If the second indicator does not change, then
Remark~\ref{comm:7} shows that $\delta'(B) \le 0$ and,
hence, the third indicator does not increase.

If the second and third indicators do not change, then Remark~\ref{comm:7} shows that
$A' = A^- \subset A$, and, hence, the fourth indicator decreases by at least~1.

Hence, $(A,B) \succ (A',B')$.

If $(A',B') \in \mathcal{S}_2$ then,
by the induction hypothesis, {\normalfont Assumptions~\ref{ass:key},}  $(A',B')$ is sound,
and we have proved Consequences~\ref{conss:3}.
\end{proof}

This completes the list of consequences. \hfill\qed
\end{hypotheses}

A substantial part of the proof of the following result is
similar to the proof of Theorem~5 in~\cite{Kemperman56}.

\begin{lemma}\label{lem:e}  With {\normalfont Assumptions~\ref{ass:key},~\ref{ass:ge},~{\it and}~\ref{ass:Q1cols}},
if $\blocks(A{\cdot_2}B)\! =\!\blocks(AB),$ then $(A,B)$ is sound.
\end{lemma}

\begin{proof}
 Consider the possibility that, for all $b_1$, $b_2$ in $B$, we have $Ab_1 = Ab_2$.
Let $$L:= \gen{b_1 \overline b_2 \mid b_1, b_2 \in B} \le G.$$
Here, $AL = A$.  Consider any $(a,b) \in A \times B$.
Then, $AB \supseteq Ab = ALb \supseteq aLb$,\linebreak and $L$ is finite.
Also, $L  \supseteq  \gen{B \overline b} \supseteq B \overline b$, and,
$\abs{L} \ge \abs{B\overline b} = \abs{B} \ge 3$,
by Assumptions~\ref{ass:ge}.
By the Sylow theorems, $L$ contains a subgroup which has order 4 or an odd prime.
Thus, $\blocks(AB) \ge 1$.  Hence, $\blocks(A{\cdot_2}B) = \blocks(AB) \ge 1$,
and $(A,B)$ is sound.

It remains to consider the case where $\blocks(AB) = 0$ and, here, by the foregoing,
there exist $b_1$ and $b_2$ in $B$ such that $Ab_1 \ne Ab_2$.

Without loss of generality, we may replace $B$ with $B\overline b_1$.
On setting $x = b_2 \overline b_1$, we have $\{1,x\} \subseteq B$
and $A \ne Ax$, and, hence, $1 \ne x$.
Let $(A',B')$ be  the Kemperman transform of $(A,B)$ with respect to $x$, as
in Definition~\ref{def:bde}.  Now
Hypotheses~\ref{hyp} apply.

By Consequences~\ref{cons:1} and~\ref{cons:2}, we may assume that
$(A',B') \in \mathcal{S}_2$.

 By Consequences~\ref{conss:3}, $A'B' \subseteq AB$,
$\hbar(A',B') \le \hbar(A,B)$ and $(A',B')$ is sound.
Since $A'{\cdot_2}B' \subseteq A'B' \subseteq AB$, we see that
$$\blocks(A'{\cdot_2}B') \,\,\, \le \,\,\, \blocks(A'B') \,\,\,  \le  \,\,\, \blocks(AB) \,\, =\,\, 0.$$
By soundness, $\hbar(A',B') \ge -4$.
Hence, $\hbar(A,B) \ge \hbar(A',B') \ge -4$, and,
hence, $(A,B)$ is sound.
\end{proof}

A substantial part of the
proof of the following result is similar to
the proof of Theorem~3 in~\cite{Kemperman56}.

\begin{lemma}\label{lem:bd}  With {\normalfont Assumptions~\ref{ass:key},~\ref{ass:ge},~{\it and}~\ref{ass:Q1cols}},
if $\blocks(A{\cdot_2}B)\! \neq\! \blocks(AB)$, then $(A,B)$ is sound.
\end{lemma}

\begin{proof}  Here, there exists some block $C$ which is contained in $AB$
but is not contained in $A{\cdot_{2}}B$.  Hence, $C \,\, \cap \,\, A{\cdot_{[=1]}}B$ is nonempty.
Let $(a,b)$ be an element of $A\times B$ such that
$ab \in C \,\, \cap \,\, A{\cdot_{[=1]}}B$.

By replacing $(A,B,a,b,C)$ with
$(\overline a A, B \overline b, 1,1,  \overline a C \overline b)$,
we may assume that
$(a,b) = (1,1)$.  In particular,
$
1 \in A{\cdot_{[=1]}}B \,\, \cap\,\, C.
$
By Assumptions~\ref{ass:Q1cols}(ii) and~(i), $A - \{1\}$ and $B-\{1\}$ are  subsets
of $A{\cdot_2}B$.

Consider first the case where $A-\{1\}$ and $B-\{1\}$ are disjoint.
Then
\begin{align*}
 \abs{AB} + \abs{A{\cdot_2}B} \ge 2\abs{A{\cdot_2}B} \ge 2( \abs{A - \{1\}} + \abs{B-\{1\}})= 2\abs{A} + 2\abs{B} -4,
\end{align*}
and $(A,B)$ is sound.  Therefore, we may assume that
$A - \{1\}$ and $B-\{1\}$ are not disjoint and, hence, there exists some $x \in (A \cap B) - \{1\}$.

Since $1 \in  A{\cdot_{[=1]}}B$, we see that $1 \in A1  - Ax$.  In particular,  $A \ne Ax$ and $1 \ne x$.
Let $(A',B')$ be  the Kemperman transform of $(A,B)$ with respect to $x$, as
in Definition~\ref{def:bde}.  Now Hypotheses~\ref{hyp} apply.

By Consequences~\ref{cons:1} and~\ref{cons:2}, we may assume that
$(A',B') \in \mathcal{S}_2$.

 By Consequences~\ref{conss:3}, $A'B' \subseteq AB$,
$\hbar(A',B') \le \hbar(A,B)$ and $(A',B')$ is sound.
By Definitions~\ref{defs:sound}, there are three possibilities.

\medskip

\noindent\textbf{Case 1.}  $\hbar(A',B') \ge -4$.

Here, $ \hbar(A,B) \ge \hbar(A',B') \ge -4$, and $(A,B)$ is sound.

\medskip

\noindent\textbf{Case 2.}   $\blocks(A'{\cdot_2}B') \ge 1$.

Here, $A'{\cdot_2}B'$ contains some block, $D$.

We claim that $C \ne D$.
Since $1 \in C$ and $D \subseteq A'{\cdot_2}B'$,
it suffices to show that $1 \not\in A'{\cdot_2}B'$.

Notice that $1 = 1{\cdot}1 \in (A^+){\cdot}(B^-)$.

By Remark~\ref{comm:4},
$1 \,\,\,\not\in\,\,\,  (A^-){\cdot_2}(B^+) \,\, \cap\,\,\ A{\cdot_{[=1]}}B.$
Since $1 \in A{\cdot_{[=1]}}B$, we see that $1 \not\in (A^-){\cdot_2}(B^+)$.

 Similarly, $1 = 1{\cdot}1\in (A^-){\cdot}(B^+)$ and, by Remark~\ref{comm:5}, $1 \not\in (A^+){\cdot_2}(B^-)$.

Hence, $1 \not\in  (A'){\cdot_2}(B')$, as desired.

Thus, $C$ and $D$ are two different blocks which are contained in $AB$.

Hence, $\blocks(AB) \ge 2$  and $(A,B)$ is sound.

\medskip

\noindent\textbf{Case 3.}  $\blocks(A'B') \ge 2$.

Since $AB \supseteq A'B'$, we see that $\blocks(AB) \ge \blocks(A'B') \ge 2$, and $(A,B)$ is sound.
\end{proof}

By Lemmas~\ref{lem:e} and~\ref{lem:bd},  the induction argument is complete,
and we have proved the following.

\begin{theorem}
Every element $(A, B)$ of $\mathcal{S}_2$ is sound.
\end{theorem}

By Lemma~\ref{lem:sound}, we have the key inequality.

\begin{theorem}\label{th:presubmain}   Let  $A$ and $B$ be
finite subsets of a group $G$.  If $\abs{A} \ge 2$ and $\abs{B} \ge 2$, then
$\abs{AB} + \abs{A{\cdot_2}B}\,\,\ge\,\, \min\{ 2\abs{A} + 2\abs{B} - 4, \, 2{\cdot}\height(G)\}.$ \hfill\qed
\end{theorem}

The proof of Corollary~\ref{cor:inequalities}(ii) is now complete.

\section{Upper bounds}\label{sec:freeproduct}

In this section, we use the viewpoint of Mihalis Sykiotis~\cite[Proof of Theorem~2.13(1)]{Sykiotis05}
together with Corollary~\ref{cor:inequalities}(ii) to rewrite and generalize results of~\cite{Ivanov99} and~\cite{Ivanov01}.

The following is well known and easy to prove.

\begin{lemma}\label{lem:prod} Let  $H$ and $K$ be subgroups of a group $G$,
and let $S$ be a set of $(H,K)$-double coset
representatives in $G$.
Then the  map
\begin{center}
$\mathop{\textstyle\bigvee}\limits_{s \in S} ((H^s \cap K)\backslash G) \quad
\to \quad (H\backslash G) \times (K\backslash G),
\qquad (H^s \cap K)g \mapsto (Hsg,Kg),$
\end{center} is bijective.
The inverse map is given by
$(Hx,Ky) \mapsto (H^s\cap K)ky$ for the unique $s \in S$ such that $Hx\overline{y}K = HsK$,
and any $k \in K$ such that $Hx\overline{y} = Hsk$; here $(H^s\cap K)k$ is unique.
\hfill \qed
\end{lemma}

It is convenient to recall the following.

\begin{review}\label{rev:} Suppose that $H$ is a group and that
$T$ is an $H$-free $H$-tree, that is, $H$ acts freely on $T$.

Then, with respect to any basepoint, the fundamental group
of the quotient graph $H\backslash T$ is isomorphic to $H$;
see, for example,~\cite[Corollary~I.4.2]{DicksDunwoody89}.
In particular, $H$ is a free group.

The {\it core} of $H\backslash T$, denoted $\core(H\backslash T)$, is the subgraph of
$H\backslash T$ consisting of all
those vertices and edges which lie in cyclically reduced closed paths in~$H\backslash T$.

Let $X = \core(H\backslash T)$.  We write $V\!X$ and $E\!X$ for the vertex-set and edge-set of $X$,
respectively. Every vertex of $X$ has valence at least two.

If $H$ is trivial, then $H\backslash T$ is the tree $T$, and $X$ is empty.

Now suppose that $H$ is nontrivial.

Then $H\backslash T$ is not a tree, and $X$ is nonempty and
its fundamental group is  isomorphic to $H$.
Moreover, $H$ is finitely generated if and only if $X$ is finite.

Suppose further that $H$ is finitely generated, or, equivalently, that $X$ is finite.

For each $v \in V\!X$, let $\deg_X(v)$ denote the valence of $v$ in $X$. Then
\begin{align*}
\textstyle \sum\limits_{v \in V\!X} (\deg_X(v) -2) &=
\textstyle(\textstyle\sum\limits_{v \in V\!X} \deg_X(v)) - (\textstyle\sum\limits_{v \in V\!X}2)=
(\textstyle\sum\limits_{e \in E\!X} 2) -(\textstyle\sum\limits_{v \in V\!X}2)\\
&= 2{\cdot}\abs{E\!X} - 2{\cdot}\abs{V\!X} = -2{\cdot}\chi(X) = -2{\cdot}\chi(H) = 2{\cdot}\barr(H).
\end{align*}
Thus $\barr(H) = \frac{1}{2} \textstyle\sum\limits_{v \in V\!X} (\deg_X(v) -2)$.
\hfill\qed
\end{review}

We now come to our main upper-bound result. Recall from Definitions~\ref{defs:begin}
that $\height(G)$ is $\infty$ or $4$ or an odd prime, and that
$\fun\!\height(G) = \frac{\height(G)}{\height(G)-2} \in [1,3]$.

\begin{theorem}\label{th:5}  Suppose that {\normalfont Notation~\ref{not:main}} holds.
Let $H$ and $K$ be elements of~$\mathcal{F}$,
and let $S$ be a set of    $(H,K)$-double coset
representatives in $G$.  Then
\vskip -0.8cm\null
\begin{center}
$\textstyle \sum\limits_{s \in S} \barr(H^s \cap K)
\quad \le \quad 2{\cdot}\!\fun\!\height(G){\cdot}\barr (H){\cdot}\barr(K).$
\end{center}
\end{theorem}

\begin{proof}
Clearly, we may assume that $H$ and $K$ are nontrivial.

Let $\{x_j \mid j \in J\}$ be a free generating set of $F$.

We view $G$ as the fundamental group of the following graph of groups.

Let $V = \{v_i \mid i \in I \vee \{0\}\}$, a set indexed by the disjoint union $I \vee \{0\}$.

Let $E = \{e_i \mid i \in I \vee J\}$, a set indexed by the disjoint union $I \vee J$.

Let $Z = (Z, V, E, \overline \iota, \overline \tau)$
 denote the (oriented) graph with vertex set $V$, edge  set~$E$, and  incidence relations
such that, for each $i \in I$ and $j \in J$, we have
$\overline \iota (e_i) = \overline \iota (e_j) =  \overline \tau (e_j) = v_0$,
and $\overline \tau (e_i) = v_i$.

Let $Z_0 = Z - \{e_j \mid j \in J\}$,  the unique maximal subtree of $Z$.

Let $(G(-), Z)$ be the unique graph of groups such that $ G(v_0) = \{1\}$, and
 for each $i \in I$, $G(v_i) = G_i$,  and,  for each $i \in I \vee J$,
 $G(e_i) =  \{1\}$.

In a natural way, the fundamental group $\pi(G(-), Z,Z_0)$
can be identified with the free product $F \ast \mathop{\text{\Large $*$}}\limits_{i \in I} G_i = G$.

Let $T = T(G(-),Z,Z_0)$ be the Bass--Serre tree for $(G(-),Z,
Z_0)$.  Thus $T = (T, VT, ET, \iota, \tau)$ is the $G$-graph
described as follows.

The vertex set is
$VT = \mathop{\textstyle\bigvee}\limits_{i \in I \vee \{0\}} Gv_i$, where, for each $i \in I \vee \{0\}$,
 the stabilizer $G_{v_i}$ is $G(v_i)$.

The edge set is
$ET = \mathop{\bigvee}\limits_{i \in I \vee J} Ge_i$, where, for each $i \in I  \vee J$, the stabilizer
$G_{e_i}$ is $G(e_i) = \{1\}$.

The incidence relations are such that, for each $g \in G$, $i \in I$, and $j \in J$, we have
$\iota(ge_i) = \iota(ge_j) =  gv_0$, $\tau(ge_j) = gx_jv_0$, and, $\tau(ge_i) = gv_i$.

By Bass--Serre theory, $T$ is a $G$-tree; see, for
example,~\cite[Theorem I.7.6]{DicksDunwoody89}.

Here, $G$ acts freely on the edge set $ET$, and $H$ and $K$ act freely on all of $T$.

We now use the argument in the proof of~\cite[Theorem~2.13(1)] {Sykiotis05}; see also\linebreak \cite[p.380]{Dicks94}.

We identify $ G\backslash T = Z$.

 The pullback of the two graph maps $H \backslash T \to   Z$ and
$K \backslash T \to  Z$ will be denoted
$(H \backslash T) \times_{Z} (K \backslash T)$.
As a set, $(H \backslash T) \times_{Z} (K \backslash T)$
 is a subset of\linebreak  $(H \backslash T) \times (K \backslash T)$; moreover,
$(H \backslash T) \times_{Z} (K \backslash T)$
has a natural graph structure.

We consider the map
\vskip -0.7cm\null
\begin{center}$\Phi\colon\mathop{\bigvee}\limits_{s \in S} ((H^s \cap K) \backslash T) \quad \to
\quad (H \backslash T) \times_{Z} (K \backslash T), \qquad (H^s \cap K)t \mapsto (Hst,Kt).$
\end{center}
\vskip -0.2cm
Here, $\Phi$ is a graph map.  By Lemma~\ref{lem:prod}, $\Phi$ is bijective on the
edge sets, and on the sets of vertices that map to $v_0$ in $Z$, since
$G$ acts freely on $ET \vee Gv_0$.
In particular, $\Phi$ is surjective.

Let us write
\vskip -0.7cm\null
\begin{center}$X = \core(H \backslash T)$, $Y = \core(K \backslash T)$
and $W = \mathop{\bigvee}\limits_{s \in S} \core((H^s \cap K) \backslash T)$.\end{center}
\vskip -0.3cm
Since $\Phi$ carries cores to cores, $\Phi$ induces a graph map
$\phi \colon W \to X \times_{Z} Y$.  Here, $\phi$
is injective on the edge sets, and on the sets of vertices which map to $v_0$ in $Z$.

By Review~\ref{rev:}, $X$ and $Y$ are finite and
$$\barr(H) = \textstyle\frac{1}{2} \sum\limits_{x \in V\!X} (\deg_X(x)-2), \quad
\barr(K) = \frac{1}{2} \sum\limits_{y \in V\!Y} (\deg_Y(y)-2).$$

Since $\phi$ embeds $EW$   in the finite set $EX \times_{EZ} EY$, we see that $W$ is finite, and, by Review~\ref{rev:},
$$\textstyle\sum\limits_{s \in S} \barr(H^s \cap K) = \frac{1}{2} \sum\limits_{w \in V\!W} (\deg_W(w)-2).$$

At this stage, we leave the proof of~\cite[Theorem~2.13(1)]{Sykiotis05} and
switch to the proof of~\cite[Theorem 2]{Ivanov99}.

Notice that the result we want to prove can be reformulated as
$$\textstyle \frac{1}{2}{\cdot}\sum\limits_{w \in V\!W} (\deg_W(w)-2)\le
 2{\cdot}\!\fun\!\height(G) \cdot (\frac{1}{2}{\cdot}\sum\limits_{x \in V\!X} (\deg_X(x)-2))\cdot
(\frac{1}{2}{\cdot}\sum\limits_{y \in V\!Y} (\deg_Y(y)-2)),$$
that is,
$$\textstyle \sum\limits_{w \in V\!W} (\deg_W(w)-2) \quad \le \quad \fun\!\height(G) \cdot
\sum\limits_{(x,y) \in V\!X \times V\!Y} ((\deg_X(x)-2){\cdot}(\deg_Y(y)-2)).$$

Consider any $(x,y) \in V\!X \times_{V\!Z} V\!Y$, and let
$\phi^{-1}(x,y)$ denote the preimage in $V\!W$ of $(x,y)$
under the map $\phi: V\!W \to  V\!X \times_{V\!Z} V\!Y$.
To prove the desired result, it then suffices to show that
\begin{equation}\label{eq:degs}
\textstyle \sum\limits_{w \in \phi^{-1}(x,y)} (\deg_W(w)-2) \quad \le \quad
\fun\!\height(G){\cdot}(\deg_X(x)-2){\cdot}(\deg_Y(y)-2).
\end{equation}

Let $z$ denote the common image of $x$ and $y$ in $Z$.  Thus,
there exists a unique $i \in I \vee \{0\}$ such that $z = v_i$.

\medskip

\noindent \textbf{Case 1.} $i = 0$.

We have seen that the graph map $\phi\colon W \to X \times_Z Y$ is injective on
the sets of vertices mapping to $v_0$ in $Z$.
Thus, here, $\phi^{-1}(x,y)$ consists of a single element, $w_0$, say.
Since~\eqref{eq:degs} is clear when all the $w$ have valence 2,
we may assume that $\deg_W(w_0) \ge 3$.
Recall that $\iota_W^{-1}\{w_0\}$, resp.  $\tau_W^{-1}\{w_0\}$, denotes the set of edges
of $W$ whose initial, resp. terminal, vertex is $w_0$.
 Then $$\abs{\iota_W^{-1}\{w_0\}} + \abs{\tau_W^{-1}\{w_0\}} = \deg_W(w_0).$$

It is not difficult to show that the induced
 map $\iota_W^{-1}\{w_0\} \to EZ$ is injective, and, hence,
$\iota_W^{-1}\{w_0\} \to \iota_X^{-1}\{x\}$ is injective, and, hence,
$\abs{\iota_W^{-1}\{w_0\}} \le \abs{\iota_X^{-1}\{x\}}$.
Similarly, $\abs{\tau_W^{-1}\{w_0\}} \le \abs{\tau_X^{-1}\{x\}}$.
Thus $\deg_W(w_0) \le \deg_X(x)$.

Similarly,   $\deg_Y(y) \ge \deg_W(w_0) \ge 3$.

Now we have
\begin{align*}
\textstyle \sum\limits_{w \in \phi^{-1}(x,y)} &(\deg_W(w)-2)\,\, = \,\, \deg_W(w_0) - 2
\,\, \le \,\, \deg_X(x) -2 \,\,\\&\hskip -1cm \le \,\, 1{\cdot}(\deg_X(x)-2){\cdot}(3-2)
\le \,\, \fun\!\height(G){\cdot}(\deg_X(x)-2){\cdot}(\deg_Y(y)-2),
\end{align*}
as desired.

\bigskip

\noindent \textbf{Case 2.}  $i \in I$.

Here, there exist $g_x,\,g_y \in G$   such that $x = Hg_xv_i$ and $y = Kg_yv_i$.

Notice that $\deg_X(x) = \abs{\tau^{-1}_X\{x\}}$, and that
$$\tau^{-1}_X\{x\} \subseteq  \{H\}g_xG_{i}e_i :=\{Hg_xae_i \mid a \in G_i\}.$$
Hence, there exists a subset $A$ of $G_i$ such that
$\tau^{-1}_X\{x\}  = \{H\}g_xAe_i$.  Moreover,  $A$ is unique,
since $G$ acts freely on $ET$ (on the left)  and
$G_i$ acts freely on $H \backslash G$ on the right.
Hence, $\abs{A} = \deg_X(x) \ge 2$.

Similarly, there exists a unique subset $B$ of $G_i$ such that
$\tau^{-1}_Y\{y\}  =  \{K\}g_yBe_i$, and $\abs{B} = \deg_Y(y) \ge 2$.

The embedding $\phi \colon EW \to EX \times_{EZ} EY$, gives an embedding
 $$\phi \colon\textstyle \bigvee\limits_{w \in \phi^{-1}(x,y)} \tau^{-1}_W\{w\} \,\,\,\, \into \,\,\,\,
\tau^{-1}_X\{x\} \times \tau^{-1}_Y\{y\}\,\, = \,\, \{H\}g_xAe_i \times  \{K\}g_yBe_i,$$
which, when composed with the embedding
$$\{H\}g_xAe_i \times  \{K\}g_yBe_i
 \,\,\,\,  \into  \,\,\,\,  A \times B, \quad (Hg_xae_i,Kg_ybe_i) \mapsto (a,b),$$
gives an embedding $$\textstyle\psi\colon\bigvee\limits_{w \in \phi^{-1}(x,y)} \tau^{-1}_W\{w\}
\quad \into \quad A \times B, \quad e \mapsto \psi(e).$$

Let $w \in \phi^{-1}(x,y)$.

We claim that $\psi(\tau^{-1}_W\{w\})$ is a single-quotient subset of $A \times B$, as in
Definitions~\ref{defs:var}. Let $e$, $f$ be elements of $\tau^{-1}_W\{w\}$.

There exist $s_w \in S$ and $g_w\in G$ such that
$w = (H^{s_w} \cap K)g_wv_i$.  Also, there exists a unique subset $C_w$ of $G_i$ such that
$\tau^{-1}_W\{w\} = (H^{s_w} \cap K)g_wC_we_i$, and, here,
 $\abs{C_w} = \abs{\tau^{-1}_W\{w\}} = \deg_W(w)$.
 There exist $c_e$, $c_{f}$ in $C_w$ such that
$$e = (H^{s_w} \cap K)g_wc_{e}e_i, \quad f = (H^{s_w} \cap K)g_wc_{f}e_i.$$

 Let $(a_e, b_e) = \psi(e)$, $(a_{f},  b_{f}) = \psi(f)$.  This means that, on
applying the map $\phi\colon EW \to EX \times_{EZ} EY$, we have
\begin{align*}
(Hs_wg_wc_{e}e_i, Kg_wc_{e}e_i)  &= \phi(e) = (Hg_xa_ee_i, Kg_yb_ee_i),  \\
(Hs_wg_wc_{f}e_i, Kg_wc_{f}e_i)  &=\phi(f) =  (Hg_xa_{f}e_i, Kg_yb_{f}e_i).
\end{align*}
Since $G$ acts freely on $ET$, we have
 $$(Hs_wg_wc_{e}, Kg_wc_{e})  = (Hg_xa_e, Kg_yb_e), \quad (Hs_wg_wc_{f}, Kg_wc_{f})  = (Hg_xa_{f}, Kg_yb_{f}).$$
Hence $Hg_xa_e\overline{c}_{e} = Hs_wg_w = Hg_xa_{f}\overline{c}_{f}$ and
$Kg_yb_e\overline{c}_{e} = Kg_w = Kg_yb_{f}\overline{c}_{f}$.
Since $G_i$ acts freely on the right on both $H \backslash G$ and $K \backslash G$,
we see that $a_e \overline{c}_{e} = a_{f}\overline{c}_{f}$ and
$b_e\overline{c}_{e} = b_{f}\overline{c}_{f}$.
Hence, $a_e\overline {b}_e = a_{f}\overline{b}_{f}$.

This completes the proof that $\psi(\tau^{-1}_W\{w\})$ is a single-quotient subset of $A \times B$.

Now
\begin{align*}
\textstyle \sum\limits_{w \in \phi^{-1}(x,y)} (\deg_W(w) - 2) \quad &= \quad \textstyle
\sum\limits_{w \in \phi^{-1}(x,y)} (\abs{\tau^{-1}_W\{w\}} - 2)\\& = \quad
\textstyle \sum\limits_{w \in \phi^{-1}(x,y)} (\abs{\psi(\tau^{-1}_W\{w\})} - 2) \\
&\le\quad \fun\!\height(G){\cdot}(\abs{A} -2){\cdot}(\abs{B} -2) \quad \text{by Corollary~\ref{cor:inequalities}(ii)}\\
&=\quad \fun\!\height(G){\cdot}(\deg_X(x) -2){\cdot}(\deg_Y(y) -2). \qedhere
\end{align*}
\end{proof}

For emphasis, we mention the extreme cases.

\begin{corollary}\label{cor:6}
Suppose that {\normalfont Notation~\ref{not:main}} holds.  Let $H$ and $K$ be elements of~$\mathcal{F}$,
and let $S$ be a set of    $(H,K)$-double coset
representatives in $G$.  Then the following hold.
\begin{enumerate}[\rm (i).]
\vskip-0.6cm \null
  \item $\sum\limits_{s \in S} \barr(H^s\cap K) \,\, \le \,\, 6{\cdot}\barr (H){\cdot}\barr(K).$
\vskip-0.3cm \null
  \item If $G$ is torsion-free, or, more generally, every finite subgroup of $G$ has
order at most two, then $\textstyle \sum\limits_{s \in S} \barr(H^s\cap K)\,\, \le \,\, 2{\cdot}\barr (H){\cdot}\barr(K).$
\hfill\qed
\end{enumerate}
\end{corollary}

We remark  that Corollary~\ref{cor:6}(i)  generalizes~\cite[Theorem 2]{Ivanov99}, while
Corollary~\ref{cor:6}(ii)  generalizes~\cite[Theorem 2]{Ivanov01}.

By combining Proposition~\ref{prop:lower} and Theorem~\ref{th:5}, we get our main result.

\begin{theorem}\label{th:6}
If {\normalfont Notation~\ref{not:main}} holds, then $\mathcal{F}$ is closed under
taking finite intersections.
Moreover,
$\upp(\mathcal{F})  \in [\depth(G){\cdot}\!\fun\!\height(G),\,\,2 {\cdot}\!\fun\!\height(G)]$, that is,
$$\left \{
\begin{array}{l}
\upp(\mathcal{F}) = \depth(G){\cdot}\!\fun\!\height(G) = 2 {\cdot}\!\fun\!\height(G) \text{ if $G$ has $2$-torsion; and,}\\ [-1cm] \\ \\ \\[-0.3cm]
\upp(\mathcal{F})  \in [\fun\!\height(G),\,\,2 {\cdot}\!\fun\!\height(G)] \text{ if $G$ is $2$-torsion free.}
\end{array}
\right .$$
\vskip -0.8cm \hfill\qed
\end{theorem}

We conclude by mentioning a more general problem.

\begin{remarks}\label{rems:concl}
Suppose that $G$ is a group and that $T$ is a $G$-tree.

Let $\mathcal{F}$ denote the set of those finitely generated (free)
subgroups $H$ of $G$ which have the property that, via the  restriction of the $G$-action,
$H$ acts freely on $T$.

Let $\upp(\mathcal{F})$ be defined as in~\eqref{eq:kappa}.

\medskip

(i). B.\ Baumslag~\cite{Baumslag66} showed that
if the $G$-stabilizers of the elements of $ET$ are all trivial, and
the $G$-stabilizers of the elements of $VT$ are all Howson, then
$G$ itself is Howson; equivalently,
the free product of a family of Howson groups is Howson.
Recall that $G$ is said to be  {\it  Howson} if the set of
finitely generated subgroups of $G$ is closed under finite intersections.

\medskip

(ii). It follows from Theorem~\ref{th:6} that, if the $G$-stabilizers of the elements of $ET$ are all trivial, then
 $\mathcal{F}$
is closed under finite intersections.  (The proof of Baumslag's result given
in~\cite[Theorem~1]{Ivanov99} shows this under the additional hypothesis that $G\backslash T$ is a tree.)
Here we conjectured that $\upp(\mathcal{F}) = \depth(G){\cdot}\!\fun\!\height(G)$,
and Theorem~\ref{th:6} implies that
$\upp(\mathcal{F})  \in [\depth(G){\cdot}\!\fun\!\height(G),\,\,2 {\cdot}\!\fun\!\height(G)]$.

\medskip

(iii).
D.\ E.\ Cohen~\cite[Theorem~7]{Cohen74}, generalizing Baumslag's result, showed that
if the  $G$-stabilizers of the elements of $ET$ are all finite, and
the $G$-stabilizers of the elements of $VT$ are all Howson, then
$G$ itself is Howson.

\medskip

(iv).
The proof of Cohen's result given by Sykiotis in~\cite[Corollary 2.14]{Sykiotis05}
shows that if the  $G$-stabilizers of the elements of $ET$ are all finite, then $\mathcal{F}$ is closed
under finite intersections.
(We recalled almost all of Sykiotis' argument in the above proof of Theorem~\ref{th:5}.)
Here we conjecture that   $\upp(\mathcal{F})$ is (again) given by the value
in~\eqref{eq:chars}, but our techniques
shed no light on this case.
\hfill\qed
\end{remarks}

\bigskip

\noindent{\textbf{\Large{Acknowledgments}}}

\medskip
\footnotesize

The research of the first-named author  was jointly
funded by the MEC (Spain) and the EFRD~(EU) through Projects BFM2003-06613 and MTM2006-13544.

The research of the second-named author  was partially funded
by the NSF (USA) through grants DMS 00-99612 and DMS 04-00746.

\vskip -.5cm\null

\bibliographystyle{amsplain}

\medskip

\textsc{
Warren Dicks,
Departament de Matem\`atiques,
Universitat Aut\`onoma de Barcelona,
E-08193 Bellaterra (Barcelona), Spain
}

\emph{E-mail address}{:\;\;}\texttt{dicks@mat.uab.cat}

\emph{URL}{:\;\;}\texttt{http://mat.uab.cat/$\sim$dicks/}

\medskip

\textsc{S.\ V.\ Ivanov, Department of Mathematics, University
of Illinois at Urbana\d1Cham\-paign, 1409 West Green Street,
Urbana IL 61801, USA }

\emph{E-mail address}{:\;\;}\texttt{ivanov@math.uiuc.edu}
\end{document}